\documentclass[12pt,final]{article}   

\usepackage{amsmath}
\usepackage{amsthm}
\usepackage{amssymb}
\usepackage{amsfonts}
\usepackage{latexsym}               
\usepackage{amsgen}
\usepackage[mathscr]{eucal}
\usepackage{mathrsfs}
\usepackage{enumerate}
\usepackage{bm}
\usepackage{here}

\usepackage{mathptmx} 
\usepackage{color}
\usepackage{fancybox}
\usepackage{graphicx}
\usepackage{cite}

\newtheorem{lem}{Lemma}[section]
\newtheorem{thm}[lem]{Theorem}
\newtheorem{pro}[lem]{Proposition}

\newtheorem{rem}[lem]{Remark}

\newcommand{\bqn}{\begin{equation}}
\newcommand{\eqn}{\end{equation}}
\newcommand{\beqx}{\begin{equation*}}
\newcommand{\eeqx}{\end{equation*}}
\newcommand{\barr}{\begin{array}}
\newcommand{\earr}{\end{array}}
\newcommand{\beqn}{\begin{eqnarray}}
\newcommand{\eeqn}{\end{eqnarray}}
\newcommand{\beqnx}{\begin{eqnarray*}}
\newcommand{\eeqnx}{\end{eqnarray*}}
\newcommand{\bmt}{\begin{multline}}
\newcommand{\emt}{\end{multline}}

\numberwithin{equation}{section}

\newcommand{\cR}{{\mathcal R}}

\newcommand{\ve}{\varepsilon}

\title{Far-field Layers of Steady State \\
for Viscous and Heat-conductive Gas}

\author{%
{\large\sc Akitaka Matsumura${}^1$}
{\normalsize and}
{\large\sc Masahiro Suzuki${}^{2}$}
\\[4pt]
}

\date{%
\normalsize
${}^1$%
Emeritus Professor, The University of  Osaka,   
\\
1-1 Machikaneyama-cho, Toyonaka, Osaka 560-0043, Japan
\\ [7pt]
${}^2$%
Department of Computer Science and Engineering, 
Nagoya Institute of Technology,
\\
Gokiso-cho, Showa-ku, Nagoya 466-8555, Japan
}

\begin{document}

\maketitle

\begin{abstract}
We study the vanishing heat conductivity limit of radially symmetric
stationary solutions to the non-isentropic compressible Navier--Stokes
equations in the exterior of a ball under an inflow boundary condition.
A distinctive feature of this singular limit is the formation of a
{\it far-field layer} whose transition region moves toward spatial infinity
as the heat conductivity tends to zero. We establish the existence and
uniqueness of stationary solutions with the uniform estimates with respect to the heat
conductivity and analyze the corresponding limiting problem.
The limiting solution generally fails to attain the prescribed far-field
density and temperature separately, although its pressure at infinity
agrees with the prescribed far-field pressure. We show that this
discrepancy gives rise to a far-field layer, derive its explicit
representation, and identify its characteristic spatial scale.
Finally, we prove that the stationary solution is uniformly approximated
over the entire exterior domain by the superposition of the limiting
solution and the far-field layer.
\end{abstract}

\begin{description}

\item[{\it Keywords:}]  
compressible Navier--Stokes equations, inflow problem, radially symmetric stationary solution, vanishing heat conductivity limit

\item[{\it 2020 Mathematics Subject Classification:}]
35Q30, 
76N06, 
76M45. 

\end{description}

\section{Introduction} 
In the present paper, we study the vanishing heat conductivity limit of radially symmetric stationary solutions to the non-isentropic compressible Navier--Stokes equations in the exterior of a ball in $\mathbb{R}^n$ under an inflow boundary condition. A distinctive feature of this singular limit is the emergence of a {\it far-field layer}. In contrast to a conventional boundary layer localized near the physical boundary, this layer develops in a spatial region that moves toward infinity as the heat conductivity vanishes.

The motion of a compressible, viscous, heat-conductive, Newtonian polytropic fluid is governed by the following non-isentropic compressible Navier--Stokes equations:
\begin{eqnarray}
 \label{cns} 
 \left\{
\begin{array}{l}
\rho_t+\nabla\cdot (\rho U)=0, \\[10pt]
(\rho U)_t+\nabla\cdot(\rho U\otimes U)+\nabla p =\nabla \cdot {\cal D}, 
\\[10pt]
\left\{\rho \left(e + \frac{1}{2}|U|^2\right)\right\}_t +
\nabla\cdot \left\{\rho U \left(e + \frac{1}{2}|U|^2\right)\right\}+\nabla\cdot(pU)= \nabla\cdot ({\cal D}U -q).\\[10pt]
\end{array}
\right. 
\end{eqnarray} 
Here, $t>0$ and $x=(x_{1},x_{2},\ldots,x_{n})\in\mathbb R^{n}$ denote the time and space variables, respectively, while the unknowns $\rho=\rho(t,x)>0$, $U={}^t(u_1(t,x),\ldots,u_n(t,x))\in\mathbb R^{n}$, and $\theta=\theta(t,x)>0$ denote the mass density, fluid velocity, and absolute temperature, respectively.
Furthermore, $e:=e(\theta)=c_V\theta$ and $p:=p(\rho,\theta)=R\rho\theta$ denote the internal energy per unit mass and the pressure, respectively.
The positive constants $c_V$ and $R$ denote the specific heat at constant volume and the gas constant, respectively. We also introduce the adiabatic exponent $\gamma$ defined by
\begin{gather*}
\gamma := 1 +\frac{R}{c_{V}}.
\end{gather*}
Furthermore, ${\cal D}$ and $q$ denote the viscous stress tensor and the heat flux vector, respectively, defined by
\begin{align}
{\cal D}= \nu \{\nabla\otimes U + {}^t(\nabla\otimes U)\}+\lambda (\nabla\cdot U)I,\quad q=-\kappa \nabla\theta,
\end{align}
where $\nu$ and $\lambda$ are the shear and second viscosity coefficients, respectively, and $\kappa$ is the heat conductivity.
These coefficients are assumed to be constants satisfying
\begin{align}
\label{visc}
\nu>0, \quad 2\nu+ n\lambda \ge 0, \quad \kappa>0.
\end{align}

In this paper, we consider the boundary value problem for \eqref{cns} in the exterior domain
$\Omega:=\{x\in\mathbb{R}^{n} \, ;  \, |x|>r_0\}$, where $n\geq 3$ and $r_0>0$. 
In particular, we focus on radially symmetric stationary solutions, which have the form
\begin{align}
\label{t}
\rho(t,x)=\rho(r), \quad U(t,x)=\frac{x}{r}\,u(r),\quad \theta(t,x)=\theta(r), \quad r= |x|,
\end{align}
where $u=u(r)$ is a scalar function. By substituting (\ref{t}) into (\ref{cns}), 
we obtain the equations for the radially symmetric stationary solution $(\rho, u,  \theta)$:
\begin{equation}
\label{nrs0}
\left\{\begin{array}{l}
\dfrac{ (r^{n-1}\rho u)_r}{r^{n-1}}=0, \\[10pt]
\dfrac{ (r^{n-1}\rho u^2)_r}{r^{n-1}}+ p_r =\mu \Big(\dfrac{(r^{n-1}u)_r}{r^{n-1}}\Big)_r,
\\[10pt]
\dfrac{\left\{r^{n-1}\rho u \left(e + \frac{1}{2}u^2\right)\right\}_r}{r^{n-1}}
+\dfrac{ (r^{n-1}pu)_r}{r^{n-1}} = \kappa\dfrac{(r^{n-1}\theta_r)_r}{r^{n-1}} +
\mu \Big(\dfrac{(r^{n-1}u)_r}{r^{n-1}}\Big)_ru +\Psi,\quad r\ge r_0,
\\[10pt]
\end{array}
 \right.\,
\end{equation}
where 
\begin{gather}
\notag 
\mu := 2\nu + \lambda>0,
\\
\Psi=\Psi(r,u,u_{r}) := 2\nu u_r^2 +2\nu(n-1)\left(\dfrac{u}{r}\right)^2+\lambda\left(u_r+\dfrac{n-1}{r}u\right)^2.
\end{gather} 
Note that $\Psi$ is the dissipation function, a non-negative quadratic form under the condition (\ref{visc}), representing the rise in temperature due to viscosity. 

We impose the following boundary and far-field conditions on \eqref{nrs0}:
\begin{gather}
\left\{
\begin{array}{ll}
(\rho,u,\theta)(r_{0}) = 
(\rho_-,u_-,\theta_-) & \text{if $u_{-}>0$},
\\
(u,\theta)(r_{0}) = 
(u_-,\theta_-) & \text{if $u_{-} \leq 0$},
\end{array}
\right.
\quad
\displaystyle{\lim_{r\to\infty}}(\rho,u,\theta)(r)
=(\rho_+,u_+,\theta_+),
\label{BFC1}
\end{gather}
where $\rho_{\pm}>0, u_{\pm}$, and $\theta_{\pm}>0$ are given constants. 
The case $u_->0$ is referred to as the inflow problem; the case $u_-=0$ as the impermeable wall problem; 
and the case $u_-<0$ as the outflow problem.

By integrating the first equation in (\ref{nrs0}) with respect to $r$, we see that 
\begin{align*}
r^{n-1}\rho(r) u(r)=\ve, \quad r\ge r_0, 
\end{align*}
for some constant $\ve$. This implies that 
\begin{align*}
u_+=\lim_{r\to \infty}u(r)= \lim_{r\to \infty}\frac{\ve}{r^{n-1}\rho_+}=0.
\end{align*}
Therefore, it is necessary to assume $u_+=0$ for the solvability of the problem \eqref{nrs0} with \eqref{BFC1}.

For the impermeable wall problem with $u_-=0$, the solution can be explicitly written as 
\[
\rho(r)=\rho_+\theta_+/\theta(r), \quad
u(r)=0, \quad
\theta(r)=\theta_+
+(\theta_--\theta_+)(r/r_0)^{-(n-2)}.
\]
For more details, see \cite{HM3}. Since the solution is explicitly available and independent of $\kappa$ in this case, no singular-layer analysis is required.

For the outflow problem with $u_-<0$, the vanishing heat conductivity limit is also nontrivial. In contrast to the inflow problem, however, a boundary layer is expected to emerge near $r=r_0$ rather than in the far-field region. We leave the rigorous analysis of this problem for future work.

\begin{figure}[t]
\centering
\includegraphics[width=0.95\textwidth]{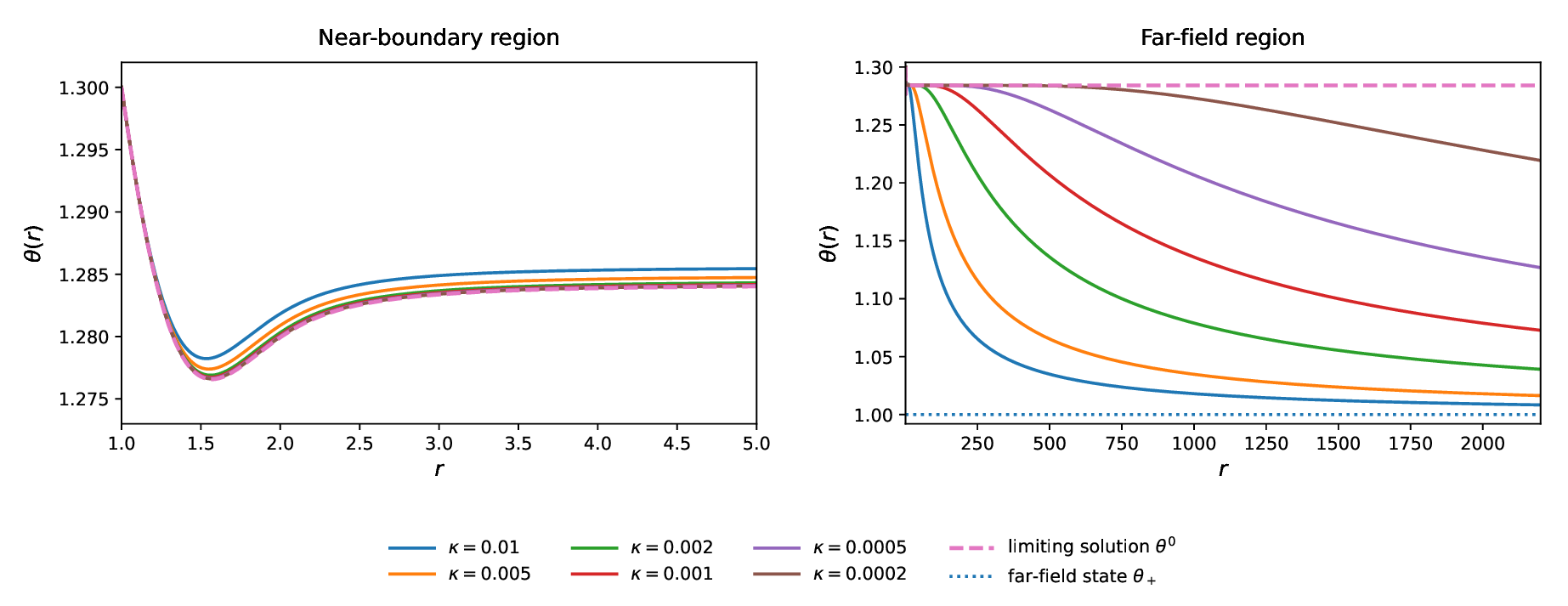}
\caption{
Temperature profiles for several values of the heat conductivity $\kappa$.
The left panel shows the behavior near the boundary, while the right panel
shows the behavior in the far-field region.
The dashed line represents the limiting solution $\theta^0$, and the dotted
line represents the prescribed far-field state $\theta_+$.
}
\label{fig:far-field-layer}
\end{figure}

For the inflow problem with $u_->0$, numerical computations clearly reveal the emergence of a far-field layer. To illustrate this behavior, we compare the stationary solutions for several values of $\kappa$ with the solution of the limiting problem obtained by formally setting $\kappa=0$, which we denote by $\theta^0$ and refer to as the limiting solution. Figure \ref{fig:far-field-layer} shows the corresponding temperature profiles.

As shown in the left panel of Figure \ref{fig:far-field-layer}, the stationary
solution approaches the limiting solution near the boundary as $\kappa$
decreases. On the other hand, the right panel shows that the stationary
solution eventually departs from the limiting solution and approaches the
prescribed far-field state $\theta_+$ as $r\to\infty$. Moreover, the transition
region moves farther away from the boundary as $\kappa$ decreases.
The numerical results therefore suggest that, as the heat conductivity tends
to zero, the stationary solution converges to the solution of the limiting
problem in the region near the boundary, whereas it approaches the constant
far-field state $(\rho_+,u_+,\theta_+)$ as $r\to\infty$.
The transition between these two asymptotic regimes occurs at a characteristic
spatial scale that moves toward infinity as the heat conductivity vanishes.
This transition region constitutes the far-field layer.

One of the main purposes of the present paper is to rigorously justify this
observation. More precisely, we identify the characteristic spatial scale at
which the transition occurs and clarify the mathematical structure of the
far-field layer. A particularly interesting feature is that the far-field
layer admits an explicit representation. Our main result then shows that, as
$\kappa\to0$, the solution of \eqref{nrs0} is uniformly approximated with
respect to $r$ by the superposition of the solution of the limiting problem
and the far-field layer. Thus, the far-field layer provides the correction
necessary to obtain a uniform approximation over the entire exterior domain.
The precise statements of our main results are given in Section \ref{S2},
while the construction and detailed analysis of the far-field layer are
carried out in Sections \ref{S6} and \ref{S7}.

We briefly review some related results on stationary solutions and
singular limits for compressible Navier--Stokes equations in exterior
domains.
For the impermeable wall problem, Matsumura and Nishida \cite{MNT}
studied the initial-boundary value problem for the non-isentropic
compressible Navier--Stokes equations in a general exterior domain in
$\mathbb{R}^{3}$. Under the boundary conditions $u_{\pm}=0$ and
$\theta_{+}=\theta_{-}$, together with
smallness assumptions on the initial perturbation, they established
the global-in-time existence of solutions around a constant state and
their large-time asymptotic behavior. For the radially symmetric
problem in dimensions $n=2,3$, Jiang \cite{J} proved the global-in-time
existence of solutions under $u_{\pm}=0$ and
$\theta_{+}=\theta_{-}$, together with partial results on the
asymptotic stability of the constant state for large initial
perturbations. Nakamura, Nishibata, and Yanagi \cite{NN} as well as Nakamura
and Nishibata \cite{NN0}, further studied related problems with
external potential forces in exterior domains and established the
existence and asymptotic stability of stationary solutions.

For the inflow and outflow problems, Hashimoto and Matsumura
\cite{HM1} studied radially symmetric stationary solutions to the isentropic
compressible Navier--Stokes equations for a viscous barotropic gas in
the exterior of a sphere and established their existence in a
sufficiently small neighborhood of the prescribed far-field state.
The asymptotic stability of these stationary solutions was subsequently
studied by Hashimoto, Nishibata, and Sugisaki \cite{HNS} for the
outflow problem and by Huang, Hashimoto, and Nishibata \cite{HHN} for
the inflow problem. For the non-isentropic heat-conductive system
considered in the present paper, the existence of radially symmetric
stationary solutions was established in \cite{HM3}.
The asymptotic stability of stationary solutions for the outflow problem was established by Huang, Nishibata, and Matsumura \cite{HNM1}.

Another closely related direction concerns singular limits of stationary solutions.
The inviscid limit of radially symmetric stationary solutions has been studied
for both the isentropic and non-isentropic compressible Navier--Stokes equations.
For the isentropic system, Hashimoto and Matsumura \cite{HM3} studied the
inviscid limit of radially symmetric stationary solutions in the exterior of
a sphere. The corresponding inviscid limit for the non-isentropic system
was studied by Hashimoto \cite{H1}. In these inviscid limits, conventional
boundary layers localized near the physical boundary $r=r_0$ arise.

The singular structure arising in the present vanishing heat conductivity
limit is fundamentally different. As $\kappa\to0$, no analogous layer is
formed near the physical boundary. Instead, the degeneration of the energy
equation causes the limiting problem to lose the prescribed far-field
conditions on the density and temperature, and the resulting discrepancy
is resolved by a far-field layer whose transition region moves toward
spatial infinity as $\kappa\to0$.
To the best of our knowledge, the formation of such a far-field layer
in the vanishing heat conductivity limit, together with its explicit
profile and a uniform composite approximation over the entire exterior
domain, has not been rigorously analyzed before.

Before closing this section, we introduce some frequently used notation.
For a nonnegative integer $k$, we denote by ${\cal B}^k(I)$
the space of functions whose derivatives up to order $k$ are continuous
and bounded on the interval $I$.
When $I=[1,\infty)$, we simply write ${\cal B}^k$ and $C^k$
for ${\cal B}^k([1,\infty))$ and $C^k([1,\infty))$, respectively.
We also denote the supremum norm on $[1,\infty)$ by
\[
\|f\|:=\sup_{r\geq1}|f(r)|.
\]
We use $c$ and $C$ to denote generic positive constants
depending on $n$, $\nu$, $\mu$, $R$, $c_V$, $\rho_+$, $\theta_+$, and $v_+$,
but independent of $r$, $\ve$, $\delta$, and $\kappa$.
We note that the positive constants $v_+$, $\ve$, and $\delta$
will be defined in the next section.
A generic positive constant depending additionally on parameters
$a,b,\ldots$ is denoted by $C_{a,b,\ldots}$.

\section{Main Results}\label{S2}
In this section, we state the main results of the present paper. Throughout this section, we consider the inflow problem $u_{-}>0$. We set $r_{0}=1$ without loss of generality and, as observed in the Introduction, necessarily have $u_{+}=0$.

We introduce the specific volume
\begin{gather}
v(r):=\frac{1}{\rho(r)},\qquad v_{\pm}:=\frac{1}{\rho_{\pm}},
\label{specific-volume}
\end{gather}
and set
\begin{gather}
\varepsilon:=\rho_{-}u_{-}=\frac{u_{-}}{v_{-}}>0,
\qquad
\delta:=\ve+|v_{-}-v_{+}|+|\theta_{-}-\theta_{+}|.
\label{epsilon-main}
\end{gather}
Here, $\ve$ represents the inflow mass flux through the boundary, while $\delta$ measures the size of the mismatch between the boundary and far-field data.

We first state the existence and uniqueness of the stationary solution for positive heat conductivity.

\begin{thm}\label{thm2.1}
Let $n\geq3$ and $\gamma<2$. For any $\rho_{+}>0$ and $\theta_{+}>0$, there exists a positive constant $\delta_{0}$, independent of $\kappa$, such that if $\delta+{\kappa}/{\varepsilon^{3}}\leq\delta_{0}$,
then the problem \eqref{nrs0} with the boundary and far-field conditions \eqref{BFC1} admits a unique radially symmetric stationary solution $(\rho,u,\theta)\in C^{1} \times C^{2} \times C^{2}$.
Moreover,
\begin{gather*}
\|v-v_{+}\|+ \|(\cdot)^{n-1}u\| +\|\theta-\theta_{+}\|
\leq C\delta,
\end{gather*}
where $C$ is a positive constant independent of $\varepsilon$, $\delta$, and $\kappa$.
\end{thm}

Although the solvability of the problem \eqref{nrs0} with \eqref{BFC1} was established in \cite{HM3}, the constants corresponding to $\delta_{0}$ and $C$ above depend on $\kappa$ in the proof given there. Since our purpose is to study the singular limit $\kappa\to0$, estimates uniform with respect to $\kappa$ are indispensable. We therefore develop a different solvability argument that provides estimates uniform in $\kappa$. The additional condition $\gamma<2$ on the adiabatic exponent is required in this argument. This condition covers the usual physically relevant polytropic gases, including monatomic, diatomic, and polyatomic gases.

We next describe the limiting state obtained by formally setting $\kappa=0$. 
Letting $\kappa=0$ in \eqref{nrs0} leads to the limiting problem
\begin{equation*}
\left\{\begin{array}{l}
\dfrac{ (r^{n-1}\rho^{0} u^{0})_r}{r^{n-1}}=0, \\[10pt]
\dfrac{ (r^{n-1}\rho^{0} (u^{0})^2)_r}{r^{n-1}}+ (R\rho^{0}\theta^{0})_r =\mu \Big(\dfrac{(r^{n-1}u^{0})_r}{r^{n-1}}\Big)_r,
\\[10pt]
\dfrac{\left\{r^{n-1}\rho^{0} u^{0} \left(c_{V}\theta^{0} + \frac{1}{2}(u^{0})^2\right)\right\}_r}{r^{n-1}}
+\dfrac{ (r^{n-1}R\rho^{0}\theta^{0} u^{0})_r}{r^{n-1}} = 
\mu \Big(\dfrac{(r^{n-1}u^{0})_r}{r^{n-1}}\Big)_ru^{0} +\Psi(r,u^{0},u^{0}_{r})
\\[10pt]
\end{array}
 \right.\,
\end{equation*}
with the conditions
\begin{gather*}
(\rho^{0},u^{0},\theta^{0})(1) = (\rho_-,u_-,\theta_-),
\quad
\lim_{r\to\infty}u^{0}(r)=u_{+}=0,
\quad
\lim_{r\to\infty}p(\rho^{0}(r),\theta^{0}(r))
=p(\rho_{+},\theta_{+}).
\end{gather*}
{ 
The far-field pressure condition is inherited from the original
momentum equation. Indeed, integrating the momentum equation
in \eqref{nrs0} from $r$ to infinity yields an identity involving
the constant $p(\rho_+,\theta_+)$, which is independent of $\kappa$.
Thus, the prescribed far-field pressure is retained in the
limiting problem, while the far-field density and temperature
are not prescribed separately.
Their individual far-field values
$(\rho^{0}_{\infty},\theta^{0}_{\infty})$
are determined by the limiting problem.
}
This loss of the prescribed far-field conditions is a consequence of the degeneration of the energy equation as $\kappa\to0$ and is precisely the mechanism that necessitates the introduction of the far-field layer.
The solvability of the limiting problem is summarized as follows.

\begin{pro}\label{pro2.2}
Let $n\geq3$. For any $\rho_{+}>0$ and $\theta_{+}>0$, there exists a positive constant $\delta_{0}$ such that if $\delta \leq\delta_{0}$, then the limiting problem admits a unique solution $(\rho^{0},u^{0},\theta^{0}) \in  C^{1} \times C^{2} \times C^{1}$. Furthermore, there exist constants $\rho^{0}_{\infty}>0$ and $\theta^{0}_{\infty}>0$ such that
\begin{gather}
 \lim_{r\to\infty}\rho^{0}(r)=\rho^{0}_{\infty},\quad
 \lim_{r\to\infty}\theta^{0}(r)=\theta^{0}_{\infty}, \quad 
p(\rho^{0}_{\infty},\theta^{0}_{\infty})=p(\rho_{+},\theta_{+}). 
\label{limiting-far-state-main}
\end{gather}
Moreover,
\begin{gather*}
 \|v^{0}-v_{+}\|+\|(\cdot)^{n-1}u^{0}\|+\|\theta^{0}-\theta_{+}\|
 \leq
 C\delta,
\end{gather*}
where $C$ is a positive constant independent of $\varepsilon$ and $\delta$.
\end{pro}

The last relation in \eqref{limiting-far-state-main} has a simple but important physical interpretation. Although the limiting solution does not, in general, approach the prescribed density and temperature $(\rho_{+},\theta_{+})$ separately, its pressure at infinity agrees exactly with the prescribed far-field pressure $p(\rho^{0}_{\infty},\theta^{0}_{\infty})=p(\rho_{+},\theta_{+})$. Thus, the pressure is already matched in the limiting problem, while the far-field layer adjusts the density and temperature so as to recover the prescribed far-field state. It is precisely this difference between $(\rho^{0}_{\infty},\theta^{0}_{\infty})$ and $(\rho_{+},\theta_{+})$ that gives rise to the far-field layer.

We now give an explicit description of this layer. Set
\begin{gather}
 v^{\rm far}(r)
 :=
 -\frac{v_{+}}{\theta_{+}}
 \bigl(\theta^{0}_{\infty}-\theta_{+}\bigr)
 \left(
 1-e^{-\frac{\varepsilon(c_{V}+R)}{\kappa(n-2)}}
 \right)^{-1}
 \left(
 e^{-\frac{\varepsilon(c_{V}+R)}{\kappa(n-2)r^{n-2}}}
 -e^{-\frac{\varepsilon(c_{V}+R)}{\kappa(n-2)}}
 \right),
\label{v-far-main}
\\
 \theta^{\rm far}(r)
 :=
 -\bigl(\theta^{0}_{\infty}-\theta_{+}\bigr)
 \left(
 1-e^{-\frac{\varepsilon(c_{V}+R)}{\kappa(n-2)}}
 \right)^{-1}
 \left(
 e^{-\frac{\varepsilon(c_{V}+R)}{\kappa(n-2)r^{n-2}}}
 -e^{-\frac{\varepsilon(c_{V}+R)}{\kappa(n-2)}}
 \right),
\label{theta-far-main}
\end{gather}
and define the corresponding velocity correction by
\begin{gather}
 u^{\rm far}(r):=\frac{\varepsilon v^{\rm far}(r)}{r^{n-1}}.
\label{u-far-main}
\end{gather}
Notice that
\begin{gather*}
 v^{\rm far}(1)=0,\quad
 \theta^{\rm far}(1)=0,
\end{gather*}
whereas
\begin{gather*}
 \lim_{r\to\infty}v^{\rm far}(r)
 =-\bigl(v^{0}_{\infty}-v_{+}\bigr),\quad
 \lim_{r\to\infty}\theta^{\rm far}(r)
 =-\bigl(\theta^{0}_{\infty}-\theta_{+}\bigr), \quad \text{where} \quad v^{0}_{\infty}:=\frac{1}{\rho^{0}_{\infty}}.
\label{far-infinity-main}
\end{gather*}
Consequently, the superposition $\bigl(v^{0}+v^{\rm far},\,u^{0}+u^{\rm far},\, \theta^{0}+\theta^{\rm far}\bigr)$
satisfies the prescribed boundary values at $r=1$ and approaches the prescribed far-field state $(v_{+},0,\theta_{+})$ as $r\to\infty$.

The explicit formulas \eqref{v-far-main} and \eqref{theta-far-main} also reveal the spatial scale of the far-field layer. The transition occurs when
\begin{gather}
 \frac{\varepsilon}{\kappa r^{n-2}}=O(1), \quad \text{that is}, \quad  r=O\left(\left(\frac{\varepsilon}{\kappa}\right)^{\frac{1}{n-2}}\right).
 \label{far-scale-main}
\end{gather}
Hence, as $\kappa\to0$, the transition region moves toward spatial infinity. This is the essential difference between the present far-field layer and a conventional boundary layer localized near the physical boundary.

We are now ready to state the main result on the vanishing heat conductivity limit.

\begin{thm}\label{thm2.3}
Let $n\geq3$ and $\gamma<2$. For any $v_{+}>0$ and $\theta_{+}>0$, there exists a positive constant $\delta_{0}$, independent of $\varepsilon$ and $\kappa$, such that if $\delta+{\kappa}/{\varepsilon^{3}} \leq\delta_{0}$,
then the stationary solution $(\rho,u,\theta)$ to the problem \eqref{nrs0} with \eqref{BFC1} satisfies
\begin{align*}
 \left\|v-v^{0}-v^{\rm far}\right\|
 +\left\|\theta-\theta^{0}-\theta^{\rm far}\right\| 
 \leq C\delta\kappa^{1/3},
 \quad
 \left\|(\cdot)^{n-1}\left(u-u^{0}-u^{\rm far}\right)\right\| \leq C\ve\delta\kappa^{1/3},
\end{align*}
where $C$ is a positive constant independent of $\varepsilon$, $\delta$ and $\kappa$.
\end{thm}

Theorem \ref{thm2.3} shows that the solution for positive heat conductivity is uniformly approximated, over the entire exterior domain, by the superposition of the limiting solution and the explicitly given far-field layer. In particular, the limiting solution alone does not provide a uniform approximation because its state at spatial infinity generally differs from the prescribed far-field state. The far-field layer exactly compensates for this discrepancy and connects the limiting far-field state $(v^{0}_{\infty},\theta^{0}_{\infty})$ to the prescribed state $(v_{+},\theta_{+})$ on the spatial scale \eqref{far-scale-main}. Moreover, the estimate is uniform for all $r\geq1$, including the transition region whose location diverges as $\kappa\to0$.

\begin{rem} \rm
The exponent $1/3$ in Theorem \ref{thm2.3} is not intended to be optimal.
In fact, there is a trade-off between the smallness assumption on
$\kappa$ relative to $\ve$ and the convergence rate in the vanishing
heat conductivity limit. By modifying the auxiliary parameter $q<1$ in the
kernel estimate \eqref{ker1}, one can weaken the condition on
$\kappa/\ve^{3}$ at the expense of a slower convergence rate, or,
conversely, obtain a faster convergence rate under a stronger smallness
condition. Since there does not appear to be a distinguished optimal
choice among these possibilities, we choose $q=1/2$ in the present paper,
which leads to the assumption and the convergence rate stated in
Theorem \ref{thm2.3}.
\end{rem}

We conclude this section by outlining the organization of the remainder
of the paper. In Section~\ref{S3}, we reformulate the stationary problem
and its limiting problem in forms suitable for the subsequent analysis,
and establish preliminary estimates for the associated kernel functions.
In Section~\ref{S4}, we establish the existence and uniqueness of
solutions to the reformulated stationary problem, which in turn gives
the solvability of the original stationary problem uniformly with
respect to the heat conductivity.
In Section~\ref{S5}, we analyze the limiting problem and characterize
the far-field state selected by the limiting equations.
In Section~\ref{S6}, we construct the far-field layer, derive its
explicit representation, and establish the estimates needed for the
singular-limit analysis.
Finally, in Section~\ref{S7}, we estimate the remainder after subtracting
the limiting solution and the far-field layer from the stationary
solution and establish the uniform composite approximation over the
entire exterior domain.

\section{Reformulation and Preliminary}\label{S3}
\subsection{Reformulation}
This section is devoted to the reformulation of  the problem \eqref{nrs0}.
We may suppose $r_0=1$ without loss of generality. 
To reformulate the problem (\ref{nrs0}), we introduce the specific volume $v$ by $v=1/\rho$ (accordingly, denote $v_\pm$ by $1/\rho_\pm$). Using the first equation of  (\ref{nrs0}), the velocity $u$ is given in terms of $v$ as
\begin{align}
\label{utov}
u(r)=\frac{\ve}{r^{n-1}}v(r),\quad r\ge 1,
\end{align}
where $\ve = u_-/v_- $.
We further introduce the new unknown functions $\eta$ and $\chi$, as the deviation of $v$ and $\theta$ from the far field state $v_+$ and $\theta_+$, respectively, by
\begin{equation}
\label{transec}
\eta(r) = v(r) - v_+,\quad \chi=\theta(r)-\theta_+, \quad r \ge 1. 
\end{equation}
Substituting (\ref{utov}) into the second equation of (\ref{nrs0}) by using (\ref{transec}), and integrating the resultant equation with respect to $r$ on $[r, \infty)$, we can rewrite the second equation of (\ref{nrs0}) as 
\begin{equation}
\label{etar}
\begin{array}{l}
\eta_r+\dfrac{R\theta_+r^{n-1}}{\ve \mu v_+^2}\eta = \dfrac{Rr^{n-1}}{\ve  \mu v_+}\chi -  \dfrac{Rr^{n-1}\chi\eta}{\ve  \mu v_+(v_++\eta)} + \dfrac{r^{n-1}}{\ve \mu}  F(r,\eta),
\end{array}
\end{equation}
where
\begin{equation}
\label{etaF}
\begin{array}{l}
F(r,\eta):=
\dfrac{R\theta_+ \eta^2}{v_+^2(v_++\eta)}
+\dfrac{\ve^{2} v_+}{2 r^{2(n-1)}}+\dfrac{\ve^{2}\eta}{ r^{2(n-1)}}
-\ve^{2}(n-1)\displaystyle{\int_r^{\infty}\dfrac{\eta(s)}{s^{2n-1}}ds}.
\end{array}
\end{equation}
For more details of the derivation of \eqref{etar}, see \cite{HM1,HM2}.
To rewrite the third equation of (\ref{nrs0}), we multiply $r^{n-1}$ by it and substitute $u=\ve(v_++\eta)/r^{n-1}$ into the resultant equation. Then we obtain the following equation: 
\begin{align}
&\left\{\ve \left(c_V\theta + \frac{1}{2}\left(\frac{\ve (v_++\eta)}{r^{n-1}}\right)^2\right)\right\}_r
+(R\ve\theta)_r
\notag \\
&= \kappa(r^{n-1}\theta_r)_r+\ve^2\mu \Big(\dfrac{(v_++\eta)\eta_r}{r^{n-1}}\Big)_r
-2\ve^2\nu (n-1)  \left( \dfrac{(v_++\eta)^{2}}{r^{n}} \right)_{r}.
\end{align}
By integrating the above equation with respect to $r$ on $[r, \infty)$ and substituting \eqref{etar} into $\eta_{r}$, we have
\begin{align*}
&\ve c_V\chi + \ve R \chi + \dfrac{\ve^3}{2}\Big(\dfrac{v_++\eta}{r^{n-1}}\Big)^2
\\
&=  \kappa r^{n-1}\chi_r  
+\ve^2\mu \dfrac{(v_++\eta)\eta_r}{r^{n-1}}
-2\ve^2\nu (n-1)   \dfrac{(v_++\eta)^{2}}{r^{n}} 
+ \alpha,
\\
&=  \kappa r^{n-1}\chi_r  
- \ve  \dfrac{R\theta_{+}(v_++\eta)\eta}{v_{+}^{2}}
+ \ve  R \chi
+\ve  (v_++\eta) F
-2\ve^2\nu (n-1)   \dfrac{(v_++\eta)^{2}}{r^{n}}
+ \alpha,
\end{align*}
where $\alpha$ denotes the constant of integration. 
Thus, we finally have the following reformulated problems of (\ref{nrs0}) in terms of $(\eta, \chi, \alpha)$:
\begin{align}
\label{reforme}
 \left\{
\begin{aligned}
  &\eta_r+\dfrac{R\theta_+r^{n-1}}{\ve \mu v_+^2}\eta = \dfrac{Rr^{n-1}}{\ve  \mu v_+}\chi  - \dfrac{Rr^{n-1}\chi\eta}{\ve  \mu v_+(v_++\eta)} +  \dfrac{r^{n-1}}{\ve \mu}  F(r,\eta),
\\
&\chi_r 
-\frac{\ve c_V}{\kappa r^{n-1}}\chi
= 
-\frac{\alpha}{\kappa r^{n-1}} 
+2\nu (n-1)   \dfrac{\ve^2 v_+^{2}}{ \kappa r^{2n-1}}
+ \dfrac{\ve R\theta_{+}}{\kappa v_{+} r^{n-1}}\eta
+  \frac{\ve}{\kappa r^{n-1}} H(r,\eta),  \quad r>1,\\
&\displaystyle{\lim_{r\to \infty}}\eta(r) =0,\quad \displaystyle{\lim_{r\to \infty}}\chi(r) =0,  \\
  &\eta(1)= \eta_-:=v_--v_+,\quad \chi(1)= \chi_-:=\theta_--\theta_+,
 \end{aligned}
 \right.
\end{align}
where
\begin{gather}
\begin{aligned}
H(r,\eta)&:= 
 \frac{\ve^2}{2}\frac{(v_++\eta)^2}{r^{2(n-1)}} 
+   \dfrac{R\theta_{+}\eta^{2}}{v_{+}^{2}}
+2\ve\nu (n-1)   \dfrac{2v_+\eta+\eta^{2}}{ r^{n}}
 - (v_++\eta) F(r,\eta)
 \\
& = \frac{\ve^2}{2}\frac{(v_++\eta)^2}{r^{2(n-1)}} 
+2\ve\nu (n-1)   \dfrac{2v_+\eta+\eta^{2}}{ r^{n}}
\\
& \quad  - (v_++\eta) 
\left(
\dfrac{\ve^{2} v_+}{2 r^{2(n-1)}}+\dfrac{\ve^{2}\eta}{ r^{2(n-1)}}
-\ve^{2}(n-1)\displaystyle{\int_r^{\infty}\dfrac{\eta(s)}{s^{2n-1}}ds}
\right).
\end{aligned}
\label{chiH}
\end{gather}
The integration constant $\alpha$ is also unknown.
Once the desired solution $(\eta, \chi,\alpha)$ of (\ref{reforme}) is obtained, 
the velocity $u$ is immediately obtained by (\ref{utov}) and (\ref{transec}) as
\begin{equation*}
u(r)=\frac{u_-(v_++\eta(r))}{v_-r^{n-1}}.
\end{equation*} 
Therefore, to prove Theorem \ref{thm2.1}, it is sufficient to establish Theorem \ref{thm3.1}, where
\begin{gather*}
\delta=\ve+|v_{-}-v_{+}|+|\theta_{-}-\theta_{+}|=\ve+|\eta_-|+|\chi_{-}|.
\end{gather*}
\begin{thm}\label{thm3.1}
Let $n \ge 3$ and $\gamma<2$. For any $v_+>0$ and $\theta_+>0$, there exists a positive constant $\delta_0$,  independent of $\kappa$, such that if $\delta + \kappa/\ve^{3} \leq \delta_0$, then the problem \eqref{reforme} admits a unique solution $(\eta,\chi,\alpha) \in C^{2} \times  C^{2} \times \mathbb R$ that satisfies $\|(\eta,\chi)\|+ |\alpha|/\ve \leq C \delta$, where $C$ is a positive constant independent of $\ve$, $\delta$, and $\kappa$.
\end{thm}

The proof of Theorem \ref{thm3.1} requires a different approach from that used in \cite{HM3}.
In \cite{HM3}, since no uniformity with respect to $\kappa$ is required, the linear term involving $\eta$ on the right hand side of the second equation in \eqref{reforme} is evaluated at the preceding iteration step. As a consequence, the two equations in the linearized problem are essentially decoupled and can be solved successively. For the present purpose, however, such a linearization does not yield estimates uniform in $\kappa$. To obtain $\kappa$-uniform estimates, the linear term involving $\eta$ must be kept at the same iteration level. The resulting linearized problem is therefore a genuinely coupled system. This requires a different solvability argument from that used in \cite{HM3}. The details are given in Section \ref{S4}.

The limiting problem associated with \eqref{reforme} is given by
\begin{align}
\label{limpro}
 \left\{
\begin{aligned}
&\eta^{0}_r+\dfrac{R\theta_+r^{n-1}}{\ve \mu v_+^2}\eta^{0} = \dfrac{Rr^{n-1}}{\ve  \mu v_+}\chi^{0}  -\dfrac{Rr^{n-1}\chi^{0}\eta^{0}}{\ve  \mu v_+(v_++\eta^{0})}+ \dfrac{r^{n-1}}{\ve \mu}  F(r,\eta^{0}),
\\
&\frac{\ve c_V}{r^{n-1}}\chi^{0}
= 
\frac{\alpha^{0}}{ r^{n-1}} 
-2\nu (n-1)   \dfrac{\ve^2 v_+^{2}}{ r^{2n-1}}
-\dfrac{\ve R\theta_{+}}{ v_{+} r^{n-1}}\eta^{0}
-  \frac{\ve}{ r^{n-1}} H(r,\eta^{0}),  \quad r>1,\\
  & \eta^{0}(1)= \eta_-,\quad \chi^{0}(1)= \chi_-,
 \end{aligned}
 \right.
\end{align}
where $F$ and $H$ are defined in \eqref{etaF} and \eqref{chiH}, respectively.
To prove Proposition \ref{pro2.2}, it suffices to establish Proposition \ref{pro3.2}, whose proof is deferred to Section \ref{S5}.

\begin{pro}\label{pro3.2}
Let $n\ge3$. For any $v_+>0$ and $\theta_+>0$, there exist positive constants
$\delta_0$ such that if $\delta\le\delta_0$, then the limiting problem
\eqref{limpro} admits a unique solution $(\eta^0,\chi^0,\alpha^{0})\in C^{2} \times  C^{2}  \times \mathbb R$ that
satisfies $\|(\eta^0,\chi^0)\|+ |\alpha^{0}|/\ve \le C\delta$, 
where $C$ is a positive constant independent of $\ve$ and $\delta$.
Moreover, the far-field limits 
\begin{gather}\label{eta-chi-infty1}
\eta^{0}_\infty=\lim_{r\to\infty}\eta^0(r),  \quad \chi^{0}_\infty=\lim_{r\to\infty}\chi^0(r)
\end{gather}
exist and satisfy 
\begin{gather}
\theta_+\eta_{\infty}^{0}
=
v_+\chi_{\infty}^{0}.
\label{farfield-relation2}
\end{gather}
The following weighted derivative estimate also holds:
\begin{gather}\label{Deta0chi0}
\|(r^{n+1}\eta_r^0,r^{n+1}\chi_r^0)\|
\le
C\frac{\delta}{\ve}.
\end{gather}
\end{pro}

The far-field layer associated with the problem \eqref{reforme} is explicitly given as
\begin{gather}\label{FFL1}
\begin{aligned}
\eta^{\rm far}&:=
-\frac{v_+}{\theta_+}
\chi^{0}_\infty
\left(
1-e^{-\frac{\varepsilon\sigma^{\rm far}}
{\kappa(n-2)}}
\right)^{-1}
\left(
e^{-\frac{\varepsilon\sigma^{\rm far}}
{\kappa(n-2)r^{n-2}}}
-
e^{-\frac{\varepsilon\sigma^{\rm far}}
{\kappa(n-2)}}
\right),
\\
\chi^{\rm far}
&
:=
-\chi^{0}_\infty
\left(
1-e^{-\frac{\varepsilon\sigma^{\rm far}}
{\kappa(n-2)}}
\right)^{-1}
\left(
e^{-\frac{\varepsilon\sigma^{\rm far}}
{\kappa(n-2)r^{n-2}}}
-
e^{-\frac{\varepsilon\sigma^{\rm far}}
{\kappa(n-2)}}
\right),
\\
\alpha^{\rm far}
&:=
\varepsilon(c_V+R)\chi^{0}_\infty
e^{-\frac{\varepsilon\sigma^{\rm far}}
{\kappa(n-2)}}
\left(1-
e^{-\frac{\varepsilon\sigma^{\rm far}}
{\kappa(n-2)}}
\right)^{-1},
\end{aligned}
\end{gather}
where $\sigma^{\mathrm{far}}:=c_V+R$.
The construction of the far-field layer is given in Section \ref{S6}.
To prove Theorem \ref{thm2.3} on the vanishing heat conductivity limit,
it is enough to establish Theorem \ref{thm3.3}. 

\begin{thm}\label{thm3.3}
Let $n \ge 3$ and $\gamma<2$. For any $v_+>0$ and $\theta_+>0$, there exists a positive constant $\delta_0$, independent of $\ve$, $\delta$, and $\kappa$, such that if $\delta + \kappa/\ve^{3} \leq \delta_0$, then
\begin{align}
\|(\eta-\eta^{0}-\eta^{\rm far},\chi-\chi^{0}-\chi^{\rm far})\| + |\alpha-\alpha^{0}-\alpha^{\rm far}|/\ve \leq C \delta  \kappa^{1/3},
\end{align}
where $C$ is a positive constant independent of $\ve$, $\delta$, and $\kappa$.
\end{thm}
The proof of Theorem \ref{thm3.3} is deferred to Section \ref{S7}.

\subsection{Preliminary} 

Treating the right hand sides of the first two equations in \eqref{reforme} as source terms, we obtain the following integral representations of $\eta$ and $\chi$:
\begin{align}
\eta(r)
&=
G(r,1)\eta_-
+
\frac{1}{\ve\mu}
\int_{1}^{r}
G(r,s)s^{n-1} \left( \dfrac{R}{v_{+}} \chi(s) - \dfrac{R\chi(s)\eta(s)}{v_+(v_++\eta(s))} + F(s) \right)\,ds,
\label{etaF1}\\
\chi(r)
&=
 \frac{\alpha}{\kappa}
\int_{r}^{\infty}
K(r,s)\frac{1}{s^{n-1}}\,ds
\notag \\
&\quad -
 \frac{\ve}{\kappa}
\int_{r}^{\infty}
K(r,s)\frac{1}{s^{n-1}}
\left(
\frac{R\theta_{+}}{v_{+}}\eta(s)
+2\nu(n-1)\frac{\ve v_{+}^{2}}{s^{n}}
+H(s)
\right)\,ds,
\label{chiF1}
\end{align}
where the kernel functions $K$ and $G$ are defined by 
\begin{align}
{G}(r,s) &:= e^{-\frac{\omega}{\mu\ve n}(r^{n}-s^{n})},\quad {\omega}:=\frac{R\theta_+}{v_+^2},
\label{G1}\\
K(r,s) &:= e^{-\frac{\ve c_{V}}{\kappa(n-2)}(\frac{1}{r^{n-2}}-\frac{1}{s^{n-2}})}.
\label{K1}
\end{align}

For the subsequent analysis, we collect several estimates for the kernel functions arising in the integral representations of both the reformulated problem \eqref{reforme} and the limiting problem \eqref{limpro}. Although the corresponding kernels involve different coefficients, they have the same structural form. We therefore introduce the following general kernels.
For $n\ge3$, $\omega_*>0$, and $\sigma_*>0$, we define the kernels
\begin{align*}
G_*(r,s)
&:=
e^{-\frac{\omega_*}{n}(r^n-s^n)},
\qquad
r,s\geq 1,
\\
K_*(r,s)
&:=
 e^{
-\frac{\sigma_*}{n-2}
\left(
\frac1{r^{\,n-2}}-\frac1{s^{\,n-2}}
\right)},
\qquad
r,s\geq 1.
\end{align*}

The kernels appearing in \eqref{etaF1}--\eqref{chiF1}, as well as those arising in the analysis of the limiting problem \eqref{limpro} and the far-field layer, are obtained from $G_*$ and $K_*$ by appropriate choices of $\omega_*$ and $\sigma_*$. Thus, the estimates below provide a unified framework for the analysis of these kernels.

\begin{lem}[\cite{HM1,HM3}]
\label{hyoka}
For the kernal $G_{*}$, the following (i)--(v) hold:
\begin{enumerate}[(i)]
\item If $\frac{\omega_{*}}{(n-2)} \geq 1$, then it holds that  
\begin{align}
\label{hyoka1}
r^{n-2}G_{*}(r,1)=r^{n-2}e^{-\frac{\omega_{*}}{n}(r^n-1)}\le 1, \qquad r\ge 1.
\end{align}
\item 
There holds that
\begin{gather}
\label{hyoka2}
G_{*}(r,1) \leq \frac{4n^{2}}{e^{2}\omega_{*}^{2}} \frac{1}{(r^{n}-1)^{2}}.
\end{gather}
\item If  $l\ge -n$, 
then there exists a constant $C_0>0$ independent of $\omega_{*}$  such that for $f \in C^0$ satisfying $r^{l} f \in {\cal {B}}^0$,
\begin{align}
\label{ineqG}
\begin{aligned}
&r^{l+n-1}  \int_1^r G_{*}(r,s) | f | (s) ds 
\le \frac{C_0}{\omega_{*}} \| r^{l}  f\|,\quad r\ge 1.
\end{aligned}
\end{align}
\item For any $a \in \mathbb R $ and $f \in {\cal {B}}^0$, there holds that
\begin{gather}\label{ineqG2}
| G_{*} (r,1) a| \leq |a| , \quad \left| \int_{1}^r G_{*} (r,s) s^{n-1}f (s)  \,ds \right|
\leq \frac{1}{\omega_{*}} \|f\|.
\end{gather}
\item Let $f \in {\cal {B}}^0$. If $\lim_{r \to \infty} f(r) = a $, then
\begin{gather}\label{convG1}
\lim_{r \to \infty} \int_{1}^r G_{*} (r,s) s^{n-1}f (s)  \,ds = \frac{a}{\omega_{*}}.
\end{gather}
\end{enumerate}
\end{lem}
\begin{proof}
One can show the assertion (i), (iii), and (iv) following \cite{HM1,HM3}. 
Assertion (ii) follows immediately from $\sup_{x\ge0}x^2e^{-x}=4/e^2$.
It remains to show the assertion (iv).
Due to ${\lim_{r\to \infty}}f(r) =a$, for any $\delta>0$, there exists $r_{*}$ such that if $r \geq r_{*}$, then $|f(r)-a| \leq \delta$. Then we observe that
\begin{align*}
&\left|
\int_{1}^{r}G_{*}(r,s)s^{n-1}f(s)\,ds
-\frac{a}{\omega_{*}}
\right|
\\
&\le
\left|
\int_{1}^{r}G_{*}(r,s)s^{n-1}(f(s)-a)\,ds
\right|
+
|a|
\left|
\int_{1}^{r}G_{*}(r,s)s^{n-1}\,ds
-\frac1{\omega_{*}}
\right|
\\
&\le
\int_{1}^{r_{*}}G_{*}(r,s)s^{n-1}|f(s)-a|\,ds
+
\delta
\int_{r_{*}}^{r}G_{*}(r,s)s^{n-1}\,ds
+
\frac{|a|}{\omega_{*}}
e^{-\frac{\omega_{*}}{n}(r^{n}-1)}
\\
&\le
\int_{1}^{r_{*}}G_{*}(r,s)s^{n-1}|f(s)-a|\,ds
+
C \delta
+
\frac{|a|}{\omega_{*}}
e^{-\frac{\omega_{*}}{n}(r^{n}-1)}
\to C \delta \quad \text{as $r \to \infty$.}
\end{align*}
This implies \eqref{convG1}.  The proof is complete.
\end{proof}

\begin{lem}
\label{hyokaK}
For the kernal $K_{*}$, the following (i)--(iii) hold:
\begin{enumerate}[(i)]
\item For any function $f \in C^{1}$ satisfying $f'(r)-\sigma_{*}r^{1-n}f(r) \leq C_{*} r^{1-n} $ for some constant $C_{*}>0$, it holds that  
\begin{align}
\label{hyokaK1}
f(r) \leq K_{*}(r,1) f(1) + \frac{C_{*}}{\sigma_{*}} (K_{*}(r,1)-1), \quad r \geq 1. 
\end{align}
\item For any function $f \in {\cal B}^{0}$, it holds that  
\begin{gather}
\label{hyokaK2}
\sup_{r \in [1,\infty)} \int_{r}^{\infty} K_{*}(r,s) \frac{|f|(s)}{s^{n-1}} ds  \leq \frac{1}{\sigma_{*}} \|f\|.
\end{gather}
\item For { $0<q<1$}, there exists a positive constant $C$ such that for all $f \in {\cal {B}}^0 \cap C^{1}$,
\begin{gather}
\begin{aligned}
& \left|  \sigma_{*} \int_{r}^{\infty} K_{*}(r,s) \frac{f(s)}{s^{n-1}} ds -f(r) \right| 
\\
& \leq \left(1+\sigma_{*} \right) e^{-\frac{\sigma_{*}}{n-2}\left(\frac{(r+\sigma_{*}^{-q})^{n-2}-r^{n-2}}{r^{n-2}(r+\sigma_{*}^{-q})^{n-2}}\right)}\|f\| 
+ \sigma_{*}^{-q} \sup_{s \in [r,r+\sigma_{*}^{-q}] } |f_{r}(s)|,
\end{aligned}
\label{ker1} 
\end{gather}
where the right hand side converges to $0$ as $\sigma_{*} \to \infty$.
In particular, when $q=1/2$ and $r=1$, we have 
\begin{gather}
\left|  \sigma_{*} \int_{1}^{\infty} K_{*}(1,s) \frac{f(s)}{s^{n-1}} ds -f(1) \right| 
 \leq { e^{(n-1)/2}}\left(1+\sigma_{*} \right) e^{-\sigma_{*}^{1/2}}\|f\| 
+ \sigma_{*}^{-1/2} \sup_{s \in [1,1+\sigma_{*}^{-1/2}] } |f'(s)|.
\label{ker2} 
\end{gather}
\end{enumerate}
\end{lem}

\begin{proof}
We first prove \eqref{hyokaK1}.  
Let $$\mu(r):=e^{\frac{\sigma_{*}}{n-2}\frac{1}{r^{n-2}}}.$$
It is easy to see from $f'(r)-\sigma_{*}r^{1-n}f(r) \leq C_{*} r^{1-n} $ that 
$
(\mu f)'(r)
\leq
C_{*} r^{1-n}\mu(r).
$
Integrating this over $[1,r]$ and using $\mu'(s)=-\sigma_{*} s^{1-n}\mu(s)$, we obtain
\begin{align*}
\mu(r)f(r)-\mu(1)f(1)
&\leq
C_{*}\int_{1}^{r} s^{1-n}\mu(s)\,ds=\frac{C_{*}}{\sigma_{*}}(\mu(1)-\mu(r)).
\end{align*}
Dividing this by $\mu(r)$ concludes \eqref{hyokaK1}.

We next prove \eqref{hyokaK2}. It is seen that
\[
\left|\int_{r}^{\infty} K_{*}(r,s) \frac{f(s)}{s^{n-1}} ds \right|
\le
\int_{r}^{\infty}
K_{*}(r,s)
\frac{1}{s^{n-1}}
\,ds \|f\|
=\frac{1}{\sigma_{*}}
\left(1-e^{-\frac{\sigma_{*}}{n-2}\frac{1}{r^{n-2}}}\right) \|f\|
\le
\frac{1}{\sigma_{*}} \|f\|.
\]
Thus \eqref{hyokaK2} holds.

Finally we show \eqref{ker1} and \eqref{ker2}. Using $\sigma_{*} s^{1-n} K_{*}(r,s)  = -   (K_{*}(r,s))_{s}$, we observe that
\begin{align*}
&f(r) - \sigma_{*} \int_{r}^{\infty} K_{*}(r,s) \frac{f(s)}{s^{n-1}} ds 
\\
&= f(r) + \int^{r+\sigma_{*}^{-q}}_{r} (K_{*}(r,s))_{s} f(s) ds
- \sigma_{*}  \int_{r+\sigma_{*}^{-q}}^{\infty} K_{*}(r,s) \frac{f(s)}{s^{n-1}} ds
\\
&= K_{*}(r,r+\sigma_{*}^{-q})f\left(r+\sigma_{*}^{-q}\right)
- \int^{r+\sigma_{*}^{-q}}_{r} K_{*}(r,s) f_{s}(s) ds
 - \sigma_{*} \int_{r+\sigma_{*}^{-q}}^{\infty} K_{*}(r,s) \frac{f(s)}{s^{n-1}} ds
 \\
&=:I_{1}+I_{2}+I_{3}.
\end{align*}
We estimate each term on the right hand side as follows:
\begin{align*}
|I_{1}| \leq & \left|K_{*}(r,r+\sigma_{*}^{-q})  f\left(r+\sigma_{*}^{-q}\right)\right| =  e^{-\frac{\sigma_{*}}{n-2}\left(\frac{(r+\sigma_{*}^{-q})^{n-2}-r^{n-2}}{r^{n-2}(r+\sigma_{*}^{-q})^{n-2}}\right)} \|f\|, 
\\
|I_{2}| \leq & \sigma_{*}^{-q} \sup_{s \in [r,r+\sigma_{*}^{-q}] } |f_{r}(s)|,
\\
|I_{3}| \leq & \sigma_{*} e^{-\frac{\sigma_{*}}{n-2}\left(\frac{(r+\sigma_{*}^{-q})^{n-2}-r^{n-2}}{r^{n-2}(r+\sigma_{*}^{-q})^{n-2}}\right)} \int_{1}^{\infty} \frac1{s^{n-1}} ds  \|f\|
\leq \sigma_{*} e^{-\frac{\sigma_{*}}{n-2}\left(\frac{(r+\sigma_{*}^{-q})^{n-2}-r^{n-2}}{r^{n-2}(r+\sigma_{*}^{-q})^{n-2}}\right)} \|f\|.
\end{align*}
From these, we arrive at \eqref{ker1}. It is straightforward to obtain \eqref{ker2} from \eqref{ker1}.
\end{proof}

\begin{lem}\label{lem:GK-mixed}
For the kernals $G_{*}$ and  $K_{*}$, the following hold:
\begin{align}
\sup_{r \in [1,\infty)}
\left|
\int_1^r
G_*(r,s)
\partial_sK_*(s,\infty)\,ds
\right|
\le
\frac1{\omega_*}
\left(\frac{2(n-1)}{e}\right)^{\frac{2(n-1)}{n-2}}
\sigma_*^{-\frac{n}{n-2}}.
\label{GKmixed}
\end{align}
\end{lem}

\begin{proof}
It is clear that
$\partial_sK_*(s,\infty)
=
\sigma_*s^{1-n}K_*(s,\infty)$.
From this, we have
\begin{align}
\left|
\int_1^rG_*(r,s)\partial_sK_*(s,\infty)\,ds
\right|
&=
\sigma_*
\int_1^r
G_*(r,s)s^{n-1}
\Bigl[
s^{-2(n-1)}K_*(s,\infty)
\Bigr]\,ds.
\label{GKproof1}
\end{align}
Set
\[
p:=\frac{2(n-1)}{n-2},
\qquad
a:=\frac{\sigma_*}{n-2}.
\]
Introducing the variable $x= a s^{2-n}$, we obtain
$
s^{-2(n-1)}K_*(s,\infty)
=
a^{-p}x^pe^{-x}.
$
The function $x^pe^{-x}$ attains its maximum at $x=p$. 
Therefore
\begin{align}
\sup_{s \in [1,\infty)}
s^{-2(n-1)}K_*(s,\infty)
\le
\left(\frac{2(n-1)}{e}\right)^{\frac{2(n-1)}{n-2}}
\sigma_*^{-\frac{2(n-1)}{n-2}}.
\label{GKproof2}
\end{align}
On the other hand, it is easy to see that
\begin{align}
\int_1^rG_*(r,s)s^{n-1}\,ds
=
\frac1{\omega_*}
\left(
1-G_*(r,1)
\right)
\le
\frac1{\omega_*}.
\label{GKproof3}
\end{align}
Combining \eqref{GKproof1}--\eqref{GKproof3}, we conclude \eqref{GKmixed}.
\end{proof}

\section{Solvability}\label{S4}
This section is devoted to the solvability of the problem \eqref{reforme}.
We begin with the following linear problem:
\begin{align}
\label{linear1}
 \left\{
\begin{aligned}
&\eta_r+\dfrac{R\theta_+r^{n-1}}{\ve \mu v_+^2}\eta = \dfrac{Rr^{n-1}}{\ve  \mu v_+}\chi  + \frac{r^{n-1}}{\ve \mu} f^{\eta},
\\
&\chi_r -\frac{\ve c_V}{\kappa r^{n-1}}\chi
= 
-\frac{\alpha}{\kappa r^{n-1}} 
+ \dfrac{\ve R\theta_{+}}{\kappa v_{+} r^{n-1}}\eta
+  \frac{\ve}{\kappa} \frac{1}{r^{n-1}} f^{\chi},  \quad r>1,\\
&\eta(1)= \eta_-,\quad \chi(1)= \chi_-, \\
 \end{aligned}
 \right.
\end{align}
where $\alpha$ is a given constant and $f^{\eta},f^{\chi} \in {\cal {B}}^0$ are given functions.
For each fixed $\alpha\in\mathbb R$, the linear system \eqref{linear1}
admits a unique global solution on $[1,\infty)$.
By the variation-of-constants formula, $\eta$ and $\chi$ admit the representations
\begin{align}
\eta(r)
&=
G(r,1)\eta_-
+
\frac{R}{\ve\mu v_+}
\int_{1}^{r}
G(r,s)s^{n-1}\chi(s)\,ds
+
\frac1{\ve\mu}
\int_{1}^{r}
G(r,s)s^{n-1}f^{\eta}(s)\,ds,
\label{lineareta1}\\
\chi(r)
&=
K(r,1)\chi_-
-
\frac{\alpha}{\kappa}
\int_{1}^{r}
K(r,s)\frac1{s^{n-1}}\,ds
+
\frac{\ve R\theta_+}{\kappa v_+}
\int_{1}^{r}
K(r,s)\frac{\eta(s)}{s^{n-1}}\,ds
+
\frac{\ve}{\kappa}
\int_{1}^{r}
K(r,s)\frac{f^{\chi}(s)}{s^{n-1}}\,ds,
\label{linearchi1}
\end{align}
where the kernels $G$ and $K$ are defined in \eqref{G1} and \eqref{K1}.

The solution $(\eta,\chi)$ of the linear problem \eqref{linear1} can be decomposed as 
\begin{gather}\label{lineardeco1}
\eta=\eta_{1}+ {\alpha}\eta_{2}, \quad \chi=\chi_{1}+ {\alpha}\chi_{2},
\end{gather}
where $(\eta_{1},\chi_{1})$ and $(\eta_{2},\chi_{2})$ solve the following problems, respectively:
 \begin{align*}
 \left\{
\begin{aligned}
  &\eta_{1r}+\dfrac{R\theta_+r^{n-1}}{\ve \mu v_+^2}\eta_{1} = \dfrac{Rr^{n-1}}{\ve  \mu v_+}\chi_{1} +  \frac{r^{n-1}}{\ve \mu} f^{\eta},
\\
&\chi_{1r} 
-\frac{\ve c_V}{\kappa r^{n-1}}\chi_{1}
=  \dfrac{\ve R\theta_{+}}{\kappa v_{+} r^{n-1}}\eta_{1} + \frac{\ve}{\kappa} \frac{1}{r^{n-1}} f^{\chi},  \quad r>1,\\
  &\eta_{1}(1)= \eta_-,\quad \chi_{1}(1)= \chi_-
 \end{aligned}
 \right.
\end{align*}
and
\begin{align*}
 \left\{
\begin{aligned}
  &\eta_{2r}+\dfrac{R\theta_+r^{n-1}}{\ve \mu v_+^2}\eta_{2} = \dfrac{Rr^{n-1}}{\ve  \mu v_+}\chi_{2},
\\
&\chi_{2r} 
-\frac{\ve c_V}{\kappa r^{n-1}}\chi_{2}
= 
-\frac{1}{\kappa r^{n-1}} 
+ \dfrac{\ve R\theta_{+}}{\kappa v_{+} r^{n-1}}\eta_{2},  \quad r>1,\\
  &\eta_{2}(1)=0,\quad \chi_{2}(1)=0.
 \end{aligned}
 \right.
\end{align*}
Here, $(\eta_{1},\chi_{1})$ and $(\eta_{2},\chi_{2})$ are independent of $\alpha$.
For $\eta_{2}$, the following lemma holds:

\begin{lem}\label{lem-eta2}
$\eta_{2}(r) \leq 0$ for any $r \geq 1$.
\end{lem}
\begin{proof}
Suppose that there exists $r_{0}>1$ such that $\eta_{2}(r_{0})>0$.
Set 
$$
r_{1}:=\sup \left\{ r \in [1, r_{0}] \, ; \, \max_{s \in [1,r]}\eta_{2}(s) \leq \frac{v_{+}}{2\ve R\theta_{+}} \right\}.
$$
By the continuity of $\eta_{2}$ and the assumption $\eta_{2}(r_{0})>0$, we have $r_{1}>1$ and $\eta_{2}(r_{1})>0$.
Treating $\eta_{2}$ as a given function and then solving the equation for $\chi_{2}$  with the boundary condition $\chi_{2}(1)=0$, we have
\begin{gather*}
\begin{aligned}
\displaystyle  \chi_{2}(r)
&  = - \dfrac{1}{\kappa} \int_{1}^{r} K(r,s) \frac{1}{s^{n-1}} ds
 + \frac{\ve R \theta_{+}}{\kappa v_{+}}\int_{1}^{r} K(r,s) \frac{\eta_{2}(s)}{s^{n-1}} ds
 \\
 &  \leq  - \dfrac{1}{2\kappa} \int_{1}^{r} K(r,s) \frac{1}{s^{n-1}} ds
 =   \dfrac{1}{2\ve c_{V} } \left(1-e^{-\frac{ c_{V} \ve }{\kappa (n-2)}(\frac{1}{r^{n-2}}-1)} \right)
 < 0
\end{aligned}
\end{gather*}
for $r \in (1,r_{1}]$, where we have used $(c_{V} \ve /\kappa) s^{1-n} K(r,s)  = -   (K(r,s))_{s}$.
Treating $\chi_{2}$ as a given function and then solving the equation for $\eta_{2}$  with the boundary conditions $\eta_{2}(1)=0$, we have
\begin{align*}
\eta_{2} (r)&= \frac{R}{\ve  \mu v_+} \int_{1}^r G (r,s) s^{n-1} \chi_{2}(s)  \,ds <0
\end{align*}
for $r \in (1,r_{1}]$.
This contradicts $\eta_{2}(r_{1}) >0$.
\end{proof}

We now show that there exists a unique constant $\alpha \in \mathbb R$ satisfying the condition ${\lim_{r\to \infty}}\eta(r) =0$ and ${\lim_{r\to \infty}}\chi(r) =0$.

\begin{lem}\label{linearsolvability1}
Let $f^{\eta},f^{\chi} \in {\cal {B}}^0$ and $\lim_{r \to \infty}f^{\eta}(r)=0$.
Then there exists a unique constant $\alpha\in\mathbb R$ such that
the corresponding solution $(\eta,\chi)$ of \eqref{linear1} satisfies
\begin{gather*}
\displaystyle{\lim_{r\to \infty}}\eta(r) =0,\quad \displaystyle{\lim_{r\to \infty}}\chi(r) =0.
\end{gather*}
Moreover, the constant $\alpha$ is characterized by
\begin{gather}\label{linearalpha1}
\begin{aligned}
\alpha &=  \ve c_{V} \left( \chi_{-} + \frac{\ve R  \theta_{+}}{\kappa v_{+}}  \int_{1}^{\infty} K(1,s) \frac{\eta(s)}{s^{n-1}} ds + \frac{\ve}{\kappa} \int_{1}^{\infty}K(1,s)\frac{f^{\chi}(s)}{s^{n-1}}\,ds \right) 
\left(1-e^{-\frac{ c_{V} \ve }{\kappa(n-2)}}  \right)^{-1}.
\end{aligned}
\end{gather}
\end{lem}
\begin{proof}
First, we show that $\eta,\chi \in {\cal {B}}^0$ for any $\alpha$.
Multiplying the first and second equations in \eqref{linear1} by $\mu \kappa^{-1} \ve^{2} \theta_{+} r^{-2(n-1)}\eta$ and $\chi$, respectively, we have
\begin{align*}
&\frac{d}{dr} \left( \frac{ \mu \ve^{2} \theta_{+}}{2 \kappa r^{2(n-1)}}\eta^{2}+\frac{1}{2}\chi^{2} \right)
+(n-1) \frac{\mu \ve^{2} \theta_{+}}{\kappa r^{2n-1}}\eta^{2}
+ \dfrac{\ve R\theta_+^{2} }{\kappa v_+^2r^{n-1}}\eta^{2} 
-\frac{\ve c_V}{\kappa  r^{n-1}}\chi^{2}
\\
&=  -\frac{\alpha}{\kappa r^{n-1}}\chi
+2\frac{\ve R\theta_+}{\kappa v_+r^{n-1}}\eta\chi
+\frac{\ve\theta_+}{\kappa r^{n-1}}\eta f^{\eta}
+\frac{\ve}{\kappa r^{n-1}}\chi f^{\chi},
\end{align*}
which, together with $f^{\eta}, f^{\chi} \in {\cal {B}}^0$, gives
\begin{gather}
\begin{aligned}
&\frac{d}{dr} \left( \frac{ \mu \ve^{2} \theta_{+}}{2 \kappa r^{2(n-1)}}\eta^{2}+\frac{1}{2}\chi^{2} \right)
-\frac{4\ve (c_V+1+ R) }{\kappa r^{n-1}}\left( \frac{ \mu \ve^{2} \theta_{+}}{2 \kappa r^{2(n-1)}}\eta^{2}+\frac{1}{2}\chi^{2} \right)
\\
&\leq 
\frac{\ve}{\kappa}
\frac{\alpha^2}{\ve^2 r^{n-1}}
+\frac{\ve v_{+}^{2}}{\kappa R r^{n-1}} (f^{\eta})^2+\frac{\ve}{\kappa r^{n-1}} (f^{\chi})^2
\leq \frac{C_{\ve,\kappa,\alpha}}{r^{n-1}}.
\end{aligned}
\end{gather}
Applying \eqref{hyokaK1} with $\sigma_{*}=4\ve (c_V+1+ R)/\kappa$ to this inequality, we have
\begin{align}
\left( \frac{ \mu \ve^{2} \theta_{+}}{2 \kappa r^{2(n-1)}}\eta^{2}+\frac{1}{2}\chi^{2} \right)(r)
& \leq   C_{\ve,\kappa,\alpha} < +\infty.
\label{lineares2}
\end{align}
Thus, $\chi \in {\cal {B}}^0$.

To prove $\eta \in {\cal {B}}^0$, we recall \eqref{lineareta1}, i.e.
\begin{align*}
\eta(r)
&=
G(r,1)\eta_-
+
\frac{R}{\ve\mu v_+}
\int_{1}^{r}
G(r,s)s^{n-1}\chi(s)\,ds
+
\frac1{\ve\mu}
\int_{1}^{r}
G(r,s)s^{n-1}f^{\eta}(s)\,ds.
\end{align*}
From this, \eqref{lineares2}, and \eqref{ineqG2} with $G_{*}=G$, it follows that
\begin{gather}\label{apriori-eta}
\|\eta\|  \leq  |\eta_-| +C (\|\chi\| + \|f^{\eta}\|) \leq |\eta_-| +C_{\ve,\kappa,\alpha}.
\end{gather}
Consequently, $\eta \in {\cal {B}}^0$.

We next determine $\alpha$ so that the far-field conditions are satisfied.
From \eqref{linearchi1},  we have 
\begin{gather}
\label{sform-chi2}
\begin{aligned}
\displaystyle  \chi(r)
& =  e^{-\frac{ c_{V} \ve }{\kappa(n-2)}(\frac{1}{r^{n-2}}-1)}\chi_{-}
+ \dfrac{\alpha}{\ve c_{V} } \left(1-e^{-\frac{ c_{V} \ve }{\kappa(n-2)}(\frac{1}{r^{n-2}}-1)} \right)
+ \frac{\ve R \theta_{+}}{\kappa v_{+}} \int_{1}^{r} K(r,s) \frac{\eta(s)}{s^{n-1}} ds
\\
& \quad 
+ \frac{\ve}{\kappa}
\int_{1}^{r}
K(r,s)\frac{f^{\chi}(s)}{s^{n-1}}\,ds.
\end{aligned}
\end{gather}
Since $\eta,f^\chi\in{\cal B}^0$, the last two terms have finite limits as $r\to\infty$.
We note that the right hand side also depends on $\alpha$ through $\eta$.
However, owing to \eqref{lineardeco1}, we have $\eta=\eta_{1}+ {\alpha}\eta_{2}$, where $\eta_{1}$ and $\eta_{2}$ are independent of $\alpha$.
Substituting this decomposition into the above formula and letting $r \to \infty$, we see that the condition $\lim_{r\to\infty}\chi(r)=0$ is equivalent to
\begin{gather*}
\begin{aligned}
\displaystyle  0
& =  e^{\frac{ c_{V} \ve }{\kappa(n-2)}}\chi_{-}
+ \dfrac{\alpha}{\ve c_{V} }  \left(1 - e^{\frac{ c_{V} \ve }{\kappa(n-2)}} 
+ \frac{\ve^{2} Rc_{V} \theta_{+}}{\kappa v_{+}} \int_{1}^{\infty} K(\infty,s) \frac{\eta_{2}(s)}{s^{n-1}} ds \right)
+ \frac{\ve R  \theta_{+}}{\kappa v_{+}} \int_{1}^{\infty} K(\infty,s) \frac{\eta_{1}(s)}{s^{n-1}} ds
\\
& \quad 
+ \frac{\ve}{\kappa}
\int_{1}^{\infty}
K(\infty,s)\frac{f^{\chi}(s)}{s^{n-1}}\,ds,
\end{aligned}
\end{gather*}
where $K(\infty,s):=\lim_{r\to\infty}K(r,s)$.
The key observation is that Lemma \ref{lem-eta2} yields
\begin{gather*}
\left(1 - e^{\frac{ c_{V} \ve }{\kappa(n-2)}} 
+  \frac{\ve^{2} Rc_{V} \theta_{+}}{\kappa v_{+}} \int_{1}^{\infty} K(\infty,s) \frac{\eta_{2}(s)}{s^{n-1}} ds \right) \leq 1 - e^{\frac{ c_{V} \ve }{\kappa(n-2)}}  <0.
\end{gather*}
Therefore, the coefficient of $\alpha$ is strictly negative, and hence
there exists a unique $\alpha\in\mathbb R$ satisfying
$\lim_{r\to\infty}\chi(r)=0$.
Thus, the constant $\alpha$ is given by
\begin{align*} 
\alpha &= \ve c_{V}  \left( e^{\frac{ c_{V} \ve }{\kappa(n-2)}}\chi_{-}
+\frac{\ve R  \theta_{+}}{\kappa v_{+}} \int_{1}^{\infty} K(\infty,s) \frac{\eta(s)}{s^{n-1}} ds +  \frac{\ve}{\kappa} \int_{1}^{\infty} K(\infty,s)\frac{f^{\chi}(s)}{s^{n-1}}\,ds \right) 
\left(e^{\frac{ c_{V} \ve }{\kappa(n-2)}} -1 \right)^{-1}
\\
&=  \ve c_{V} \left( \chi_{-} + \frac{\ve R  \theta_{+}}{\kappa v_{+}}  \int_{1}^{\infty} K(1,s) \frac{\eta(s)}{s^{n-1}} ds + \frac{\ve}{\kappa} \int_{1}^{\infty}K(1,s)\frac{f^{\chi}(s)}{s^{n-1}}\,ds \right) 
\left(1-e^{-\frac{ c_{V} \ve }{\kappa(n-2)}}  \right)^{-1},
\end{align*}
where we have used $K(\infty,s) = e^{\frac{ c_V \ve }{\kappa(n-2)}}K(1,s)$.

Finally,  applying  \eqref{convG1} to \eqref{lineareta1} and using ${\lim_{r\to \infty}}\chi(r) =\lim_{r \to \infty}f^{\eta}(r)=0$,
we obtain  ${\lim_{r\to \infty}}\eta(r) =0$. The proof is complete.
\end{proof}

Using $\lim_{r\to\infty}\chi(r)=0$, we can also express $\chi$ as 
\begin{gather}
\label{linearchi2}
\begin{aligned}
\displaystyle  \chi(r)
 = \dfrac{\alpha}{\kappa} \int_{r}^{\infty} K(r,s) \frac{1}{s^{n-1}} ds
 -\frac{\ve R \theta_{+}}{\kappa v_{+}} \int_{r}^{\infty} K(r,s) \frac{\eta(s)}{s^{n-1}} ds
  -\dfrac{\ve}{\kappa} \int_{r}^{\infty} K(r,s) \frac{f^{\chi}(s)}{s^{n-1}} ds.
\end{aligned}
\end{gather}
We are now in a position to prove Theorem \ref{thm3.1}.

\begin{proof}[Proof of Theorem \ref{thm3.1}]
We set $(\eta_{0},\chi_{0})=(0,0)$ and recursively define the sequence $\{(\eta_{m},\chi_{m})\}_{m\ge 1}$ by solving
\begin{align}
\label{apppro1}
 \left\{
\begin{aligned}
  &(\eta_{m+1})_{r}+\dfrac{R\theta_+r^{n-1}}{\ve \mu v_+^2}\eta_{m+1} = \dfrac{Rr^{n-1}}{\ve  \mu v_+}\chi_{m+1}  + \frac{r^{n-1}}{\ve\mu} {\cal R}^{\eta}_{m},
\\
&(\chi_{m+1})_{r}
-\frac{\ve c_V}{\kappa r^{n-1}}\chi_{m+1}
= 
-\frac{\alpha_{m+1}}{\kappa r^{n-1}} 
+ \dfrac{\ve R\theta_{+}}{\kappa v_{+} r^{n-1}}\eta_{m+1}
+2\nu (n-1)   \dfrac{\ve^2 v_+^{2}}{ \kappa r^{2n-1}}
+  \frac{\ve}{\kappa r^{n-1}} {\cal R}^{\chi}_{m},\\
&\displaystyle{\lim_{r\to \infty}}\eta_{m+1}(r) =0,\quad \displaystyle{\lim_{r\to \infty}}\chi_{m+1}(r) =0,  \\
&\eta_{m+1}(1)= \eta_-, \quad  \chi_{m+1}(1)= \chi_-,\\
 \end{aligned}
 \right.
\end{align}
where $\alpha_{m+1}$ is chosen so that both
$\eta_{m+1}$ and $\chi_{m+1}$ vanish at infinity, as in
Lemma \ref{linearsolvability1}. More precisely,
\begin{align}
\alpha_{m+1} &=  \ve c_{V} \left( \chi_{-} + \frac{\ve R  \theta_{+}}{\kappa v_{+}}  \int_{1}^{\infty} K(1,s) \frac{\eta_{m+1}(s)}{s^{n-1}} ds + \frac{\ve}{\kappa} \int_{1}^{\infty}K(1,s)\frac{2\nu(n-1) \ve v_{+}^{2} s^{-n} + {\cal R}^{\chi}_{m}(s)}{s^{n-1}}\,ds  \right) 
\notag \\
& \quad \times \left(1-e^{-\frac{ c_{V} \ve }{\kappa(n-2)}}  \right)^{-1},
\label{alpha-m+1}
\end{align}
and $({\cal R}^{\eta}_{m},{\cal R}^{\chi}_{m})$ are defined by
\begin{align*}
{\cal R}^{\eta}_{m} :=  -\dfrac{R\chi_{m}\eta_{m}}{v_+(v_++\eta_{m})} + F(r,\eta_{m}),
\quad 
{\cal R}^{\chi}_{m}:= H(r,\eta_{m}).
\end{align*}
Here, $({\cal R}^{\eta}_{m},{\cal R}^{\chi}_{m}) \in {\cal {B}}^0 \times {\cal {B}}^0$ and $\lim_{r \to \infty}{\cal R}^{\eta}_{m}(r)=0$ hold, and hence Lemma \ref{linearsolvability1} ensures that $(\eta_{m+1},\chi_{m+1})$ is well-defined whenever $(\eta_{m},\chi_{m}) \in {\cal {B}}^0 \times {\cal {B}}^0$, $\lim_{r\to \infty}\eta_{m}(r) = \lim_{r\to \infty}\chi_{m}(r) =0$, and $|\eta_{m}| <v_{+}/2$.
We show that $\{(\eta_{m},\chi_{m})\}_{m\ge 0}$ is a Cauchy sequence in ${\cal {B}}^0 \times {\cal {B}}^0$ for suitably small $\delta$, and its limit $(\eta,\chi) \in C^{2} \times C^{2}$ is the desired solution of the problem \eqref{reforme}.

Let us define the invariant set
\begin{gather*}
X(M):=
\left\{
(\eta,\chi)\in{\cal B}^0\times{\cal B}^0;\
\|(\eta,\chi)\|\leq M\delta \leq v_+/2
\right\}.
\end{gather*}
We seek a constant $M>0$, independent of $\ve$, $\delta$, and $\kappa$, such that
$
\{(\eta_m,\chi_m)\}_{m\geq1}\subset X(M),
$
provided that $\delta$ is sufficiently small depending on $M$.
Clearly, $(\eta_0,\chi_0)\in X(M)$ for every $M>0$.
Suppose that $(\eta_m,\chi_m)\in X(M)$ for some $M>0$.
It is straightforward to see that
\begin{gather}
\|({\cal R}^{\eta}_{m},{\cal R}^{\chi}_{m})\| \leq C_{M} \delta^{2},
\end{gather}
where $C_{M}>0$ is a constant independent of $\ve$, $\delta$, and $\kappa$ but depending on $M$.
Using \eqref{linearchi2} and \eqref{hyokaK2} with $K_{*}=K$, we observe that
\begin{align*}
\| \chi_{m+1} \|
 & \leq \left| \dfrac{\alpha_{m+1}}{\ve c_{V}} \right| 
 +\frac{R \theta_{+}}{c_{V}v_{+}} \|\eta_{m+1}\|
 +\frac{2\nu(n-1)v_{+}^{2}}{c_{V}} \ve + C_{M} \delta^{2}
 \\
 & \leq \left| \dfrac{\alpha_{m+1}}{\ve c_{V}} \right| 
 +\frac{R}{c_{V}} \|\chi_{m+1}\|
 +\frac{R \theta_{+}}{c_{V}v_{+}} |\eta_{-}|
 +\frac{2\nu(n-1)v_{+}^{2}}{c_{V}} \ve + C_{M} \delta^{2},
\end{align*}
where we have also used \eqref{lineareta1} and \eqref{ineqG2} with $G_{*}=G$ in deriving the last inequality. 
From this and the condition $\gamma=1+R/c_{V}<2$, it holds that 
\begin{gather}\label{l-chi-es}
(2-\gamma) \| \chi_{m+1} \| 
\leq \left| \dfrac{\alpha_{m+1}}{\ve c_{V}} \right| 
+  \left( \frac{R \theta_{+}}{c_{V}v_{+}} |\eta_{-}| +\frac{2\nu(n-1)v_{+}^{2}}{c_{V}} \ve \right) 
+ C_{M} \delta^{2}.
\end{gather} 
Using \eqref{lineareta1} and \eqref{ineqG2} with $G_{*}=G$ again, we also have 
\begin{gather}\label{l-eta-es}
(2-\gamma) \| \eta_{m+1} \| 
\leq (2-\gamma) |\eta_{-}|
+  \frac{v_{+}}{\theta_{+}} \left| \dfrac{\alpha_{m+1}}{\ve c_{V}} \right| 
+   \frac{v_{+}}{\theta_{+}} \left( \frac{R \theta_{+}}{c_{V}v_{+}} |\eta_{-}| +\frac{2\nu(n-1)v_{+}^{2}}{c_{V}} \ve \right) 
+ C_{M} \delta^{2}.
\end{gather} 
Recalling \eqref{alpha-m+1}, applying \eqref{ker2} together with
the first equation in \eqref{apppro1}, and assuming
$\kappa/\ve\leq1$ and $\kappa/\ve^{3}\leq1$, we obtain
\begin{gather}\label{l-alpha-es}
\begin{aligned}
&\left| \dfrac{\alpha_{m+1}}{\ve c_{V}}  - \left( \chi_{-} + \frac{ R  \theta_{+}}{c_{V}v_{+}} \eta_{-} + 2\nu(n-1) \ve v_{+}^{2} \frac{\ve}{\kappa} \int_{1}^{\infty}\frac{K(1,s)}{s^{2n-1}}\,ds\right)  \left(1-e^{-\frac{ c_{V} \ve }{\kappa(n-2)}} \right)^{-1} \right| 
\\
& \leq \tilde{C}_{0}\left(1+\frac{1}{\ve}\right) \left(\frac{\kappa}{\ve} \right)^{1/2}  \|(\eta_{m+1},\chi_{m+1})\|
+ C_{M}\delta^{2},
\end{aligned}
\end{gather}
where $\tilde{C}_{0}>0$ is a constant independent of $\ve$, $\delta$, and $\kappa$ as well as  both $\frac{\ve}{\kappa} \int_{1}^{\infty}s^{1-2n} K(1,s)\,ds$ and $(1-e^{-{ c_{V} \ve }/{\kappa(n-2)}} )^{-1}$ are  bounded independently of $\ve$ and $\kappa$.
Combining \eqref{l-chi-es}--\eqref{l-alpha-es} and taking $\kappa / \ve^{3} \ll 1$, we can find a constant $\tilde{C}_{1}>0$, independent of $\ve$, $\delta$, and $\kappa$, such that
\begin{gather}
\|(\eta_{m+1},\chi_{m+1})\| \leq \tilde{C}_{1} \delta + C_{M} \delta^{2}.
\end{gather}
Therefore, choosing $M=2\tilde{C}_{1}$ and taking $\delta$ sufficiently small, we conclude by induction that $\{(\eta_{m},\chi_{m})\}_{m\ge 1} \subset X(M)$.

We next prove the contraction estimate
\begin{gather*}
\|(\eta_{m+2}-\eta_{m+1},\chi_{m+2}-\chi_{m+1})\|
\leq
\frac12
\|(\eta_{m+1}-\eta_m,\chi_{m+1}-\chi_m)\|,
\quad m\geq1,
\end{gather*}
provided that $\delta$ is sufficiently small.
Once this estimate is established, it follows immediately that
$\{(\eta_m,\chi_m)\}_{m\ge1}$ is a Cauchy sequence in
${\cal B}^0\times{\cal B}^0$.
To prove the estimate, set
\begin{gather*}
\bar{\eta}_{m+1}:=\eta_{m+2}-\eta_{m+1}, \quad \bar{\chi}_{m+1}:=\chi_{m+2}-\chi_{m+1}, \quad \bar{\alpha}_{m+1}:=\alpha_{m+2}-\alpha_{m+1}
\end{gather*}
and see that $(\bar{\eta}_{m+1},\bar{\chi}_{m+1},\bar{\alpha}_{m+1})$ satisfies 
\begin{align*}
 \left\{
\begin{aligned}
  &(\bar{\eta}_{m+1})_{r}+\dfrac{R\theta_+r^{n-1}}{\ve \mu v_+^2}\bar{\eta}_{m+1} = \dfrac{Rr^{n-1}}{\ve  \mu v_+}\bar{\chi}_{m+1}  + \frac{r^{n-1}}{\ve\mu} ({\cal R}^{\eta}_{m+1}-{\cal R}^{\eta}_{m}),
\\
&(\bar{\chi}_{m+1})_{r}
-\frac{\ve c_V}{\kappa r^{n-1}}\bar{\chi}_{m+1}
= 
-\frac{\bar{\alpha}_{m+1}}{\kappa r^{n-1}} 
+ \dfrac{\ve R\theta_{+}}{\kappa v_{+} r^{n-1}}\bar{\eta}_{m+1}
+  \frac{\ve}{\kappa r^{n-1}} ({\cal R}^{\chi}_{m+1}-{\cal R}^{\chi}_{m}),\\
&\displaystyle{\lim_{r\to \infty}}\bar{\eta}_{m+1}(r) =0,\quad \displaystyle{\lim_{r\to \infty}}\bar{\chi}_{m+1}(r) =0,  \\
&\bar{\eta}_{m+1}(1)= 0, \quad  \bar{\chi}_{m+1}(1)= 0.\\
 \end{aligned}
 \right.
\end{align*}
We note that the linear parts of these equations are the same as those
of \eqref{apppro1}, while the boundary data are zero.
Moreover, $\bar{\alpha}_{m+1}$ is given by
\begin{align*}
\bar{\alpha}_{m+1} &=  \ve c_{V} \left(  \frac{\ve R  \theta_{+}}{\kappa v_{+}}  \int_{1}^{\infty} K(1,s) \frac{\bar{\eta}_{m+1}(s)}{s^{n-1}} ds + \frac{\ve}{\kappa} \int_{1}^{\infty}K(1,s)\frac{{\cal R}^{\chi}_{m+1}(s)-{\cal R}^{\chi}_{m}(s)}{s^{n-1}}\,ds  \right) 
\notag \\
& \quad \times \left(1-e^{-\frac{ c_{V} \ve }{\kappa(n-2)}}  \right)^{-1},
\end{align*}
Furthermore, since $\{(\eta_{m},\chi_{m})\}_{m\ge 1} \subset X(M)$, the nonlinear terms satisfy
\begin{gather*}
\|({\cal R}^{\eta}_{m+1}-{\cal R}^{\eta}_{m},{\cal R}^{\chi}_{m+1}-{\cal R}^{\chi}_{m})\| \leq C_{M} \delta \|(\eta_{m+1}-\eta_{m},\chi_{m+1}-\chi_{m})\|.
\end{gather*}
Proceeding in exactly the same manner as in the derivation of
\eqref{l-chi-es}--\eqref{l-alpha-es}, we obtain
\begin{gather*}
(2-\gamma) \| \bar{\chi}_{m+1} \| 
\leq \left| \dfrac{\bar{\alpha}_{m+1}}{\ve c_{V}} \right| 
+ C_{M} \delta \|(\bar{\eta}_{m},\bar{\chi}_{m})\|,
\\
(2-\gamma) \| \bar{\eta}_{m+1} \| 
\leq  \frac{v_{+}}{\theta_{+}} \left| \dfrac{\bar{\alpha}_{m+1}}{\ve c_{V}} \right| 
+ C_{M} \delta \|(\bar{\eta}_{m},\bar{\chi}_{m})\|,
\\
\begin{aligned}
\left| \dfrac{\bar{\alpha}_{m+1}}{\ve c_{V}}  \right| 
\leq \tilde{C}_{0}\left(1+\frac{1}{\ve}\right) \left(\frac{\kappa}{\ve} \right)^{1/2}  \|(\bar{\eta}_{m+1},\bar{\chi}_{m+1})\|
+ C_{M} \delta \|(\bar{\eta}_{m},\bar{\chi}_{m})\|.
\end{aligned}
\end{gather*}
These inequalities yield the desired contraction estimate,
provided that $\kappa/\ve^{3}$ and $\delta$ are sufficiently small.
Consequently,
$\{(\eta_m,\chi_m)\}_{m\ge1}$ is a Cauchy sequence in
${\cal B}^0\times{\cal B}^0$ and therefore converges to some
$(\eta,\chi)\in{\cal B}^0\times{\cal B}^0$.

Finally, we show that the limit
$(\eta,\chi)$ belongs to $C^2\times C^2$
and is the desired solution of \eqref{reforme}.
Note that $\eta_{m+1}$ and $\chi_{m+1}$ admit the integral
representations simliar to \eqref{lineareta1} and \eqref{linearchi2},
respectively, in terms of the kernels $G$ and $K$.
Passing to the limit as $m\to\infty$ in the integral representations, we see that
$\eta$ and $\chi$ satisfy \eqref{etaF1} and \eqref{chiF1},
respectively.
These imply that
$(\eta,\chi) \in C^2\times C^2$
and hence $(\eta,\chi)$ satisfies \eqref{reforme}.
This completes the proof.
\end{proof}

\section{Limiting Problem}\label{S5}
In this section, we discuss the solvability of the limiting problem \eqref{limpro}.
First we see from the second equation in \eqref{limpro} that the function $\chi^{0}$ can be expressed by $\eta^{0}$ as
\begin{align}
\chi^{0} 
=\frac{1}{{c_V}} \left( 
\frac{\alpha^{0}}{\ve} 
-2\nu (n-1)   \dfrac{\ve v_+^{2}}{ r^{n}} 
-\dfrac{ R\theta_{+}}{ v_{+}}\eta^{0} 
- H(r,\eta^{0}) \right).
\label{form-chi1}
\end{align}
Furthermore, evaluating the second equation at $r=1$, we see that 
\begin{gather}\label{alpha1}
\frac{\alpha^{0}}{\ve}
=  c_{V}\chi_{-} +\dfrac{R\theta_{+}}{ v_{+}}\eta_{-} +2\nu (n-1) \ve { v_+^{2}}  
+ H(r,\eta^{0})|_{r=1}.
\end{gather}
Hence, $\alpha^{0}=\alpha^{0}(\eta_{-},\chi_{-},\eta^{0})$ is a constant depending on $\eta_{-}$, $\chi_{-}$, and $\eta^{0}$.

Plugging \eqref{form-chi1} into the first equation in \eqref{limpro}, we have the following single equation for $\eta^{0}$:
\begin{gather}
\eta^{0}_r+ \dfrac{(R+R^{2}c_{V}^{-1})\theta_+r^{n-1}}{\ve \mu v_+^2} \eta^{0}
=  \dfrac{Rr^{n-1}}{\ve  \mu c_{V}v_+}\frac{\alpha^{0}}{\ve} 
-2 \nu (n-1) \dfrac{R}{  \mu c_{V}}    \dfrac{ v_+}{r}  + \frac{r^{n-1}}{\ve  \mu } \cR^{0}(r,\eta^{0},\alpha^{0}),
\label{single-eta1}
\end{gather}
where
\begin{align*}
\cR^{0}(r,\eta,\alpha) & :=  -\dfrac{R\eta}{c_{V} v_+(v_{+}+\eta)}
\left( \frac{\alpha}{\ve} 
-2\nu (n-1)   \dfrac{\ve v_+^{2}}{ r^{n}} 
-\dfrac{ R\theta_{+}}{ v_{+}}\eta
- H(r,\eta)  \right)
\\
&\quad - \frac{R}{c_Vv_+}  H(r,\eta)
+F(r,\eta).
\end{align*}
Furthermore, treating $\cR^{0}$ as a given function and then solving the differential equation \eqref{single-eta1} for $\eta$  with the boundary conditions $\eta(1)=\eta_-$, we arrive at the integral equation
\begin{gather}
\begin{aligned}
\eta^{0}(r)&=G^{0}(r,1)\eta_-
+\dfrac{R}{\ve \mu c_{V} v_+} \frac{\alpha^{0}}{\ve}  \int_{1}^rG^{0}(r,s) s^{n-1}\,ds
- 2 \nu (n-1) \dfrac{Rv_{+}}{  \mu c_{V}}  \int_{1}^rG^{0}(r,s) s^{-1}\,ds
\\
& \quad  
+ \frac{1}{\ve  \mu} \int_{1}^rG^{0}(r,s)\cR^{0}(s,\eta^{0}(s),\alpha^{0}) s^{n-1} \,ds,
\end{aligned}
\label{single-eta2}
\end{gather}
where
\begin{gather*}
G^{0}(r,s) := e^{-\frac{\omega^{0}}{\mu\ve n}(r^n-s^{n})},\quad \omega^{0}:=\frac{(R+R^{2}c_{V}^{-1})\theta_+}{v_+^2}.
\end{gather*}
We sometimes decompose $\alpha^{0}$ as
\begin{gather*}
\alpha^{0}(\eta_{-},\chi_{-},\eta^{0})
=
\alpha^{0}_{1}(\eta_{-},\chi_{-})
+
\alpha^{0}_{2}(\eta_{-},\eta^{0}),
\end{gather*}
where
\begin{align*}
\alpha^{0}_{1}=\alpha^{0}_{1}(\eta_{-},\chi_{-})
&:=
\ve c_{V}\chi_{-}
+\frac{\ve R\theta_{+}}{v_{+}}\eta_{-}
+2\nu(n-1)\ve^{2}v_{+}^{2}
\\
&\quad
-\frac{\ve^{3}}{2}(v_{+}+\eta_{-})\eta_{-}
+2\ve^{2}\nu(n-1)(2v_{+}\eta_{-}+\eta_{-}^{2}),
\\
\alpha^{0}_{2}
=\alpha^{0}_{2}(\eta_{-},\eta^{0})
&:=
(n-1)\ve^{3}(v_{+}+\eta_{-})
\int_{1}^{\infty}
\frac{\eta^{0}(s)}{s^{2n-1}}\,ds.
\end{align*}

We are now in a position to show the solvability of the integral equation \eqref{single-eta2}.

\begin{lem}\label{thm-single-eta2}
Let $n \ge 3$. For any $v_+>0$ and $\theta_+>0$, there exist positive constants $\delta_0$ and $C$ such that if $\delta \leq \delta_0$, then the problem \eqref{single-eta2} admits a unique solution $\eta^{0} \in {\cal {B}}^0$ that satisfies $\|\eta^{0}\| \leq C\delta$. 
\end{lem}
\begin{proof}
First, we define the approximate sequence $\{\eta^{0}_{m}\}_{m\ge 0}$ by
\begin{align*}
\eta^{0}_{0}(r)&:=G^{0}(r,1)\eta_- 
+\dfrac{R}{\ve \mu c_{V} v_+} \frac{\alpha^{0}_{1}}{\ve}  \int_{1}^rG^{0}(r,s) s^{n-1}\,ds
- 2 \nu (n-1) \dfrac{Rv_{+}}{  \mu c_{V}}  \int_{1}^rG^{0}(r,s) s^{-1}\,ds,
\\
\eta^{0}_{m+1}(r)&:= \eta^{0}_{0}(r) 
+\dfrac{R}{\ve \mu c_{V} v_+} \frac{\alpha^{0}_{2}(\eta_{-},\eta^{0}_{m})}{\ve}  \int_{1}^rG^{0}(r,s) s^{n-1}\,ds
\\
& \quad + \frac{1}{\ve  \mu} \int_{1}^rG^{0}(r,s)\cR^{0}(s,\eta^{0}_{m}(s),\alpha^{0}(\eta_{-},\chi_{-},\eta^{0}_{m})) s^{n-1}\,ds.
\end{align*}
We note that $\eta_{0}^{0} \in {\cal {B}}^0$ satisfies the desired boundary condition; 
$\eta^{0}_{m+1}$ also satisfies the same condition whenever $\eta^{0}_{m} \in {\cal {B}}^0$. 
We show that $\{\eta^{0}_{m}\}_{m\ge 0}$ is a Cauchy sequence in ${\cal {B}}^0$ for suitably small $\delta$, and its limit $\eta^{0} \in {\cal {B}}^0$ is the desired solution of the integral equation \eqref{single-eta2}.

Let us define an invariant set by
\begin{gather*}
X(M):=\{\eta \in {\cal {B}}^0 ;\  \|\eta\| \leq M \delta \leq v_+/2 \}.
\end{gather*}
We then seek a constant $M >0$, independent of $\mu$, $\ve$, $\eta_{-}$, and $\chi_{-}$, such that $\{\eta^{0}_{m}\}_{m\ge 0} \subset X(M)$, assuming that $\delta$ is suitably small depending on $M$.
It is easy to see from \eqref{hyoka1}, \eqref{ineqG}, and \eqref{ineqG2} with $G_{*}=G^{0}$ and $l=0$ that 
\begin{gather*}
\| \eta^{0}_{0} \| \leq M_{1} \delta,
\end{gather*}
where $M_{1}$ is a positive constant independent of $\mu$, $\ve$, $\eta_{-}$, and $\chi_{-}$.
Now we set $M:=3M_{1}$. 
Suppose that $\eta^{0}_{m}$ belongs to $X(M)$. 
We first see that
\begin{gather*}
\left| \frac{\alpha^{0}_{2}(\eta_{-},\eta^{0}_{m})}{\ve} \right| \leq C_{0} \delta^{2},
\end{gather*}
where $C_{0}$ is a positive constant depending on $M$ but independent of $\mu$, $\ve$, $\eta_{-}$, and $\chi_{-}$.
Hence,
\begin{gather*}
\left\| \dfrac{R}{\ve \mu c_{V} v_+} \frac{\alpha^{0}_{2}(\eta_{-},\eta^{0}_{m})}{\ve}  \int_{1}^{\cdot}G^{0}(\cdot,s) s^{n-1}\,ds \right\| \leq C_{1} \delta^{2},
\end{gather*}
where $C_{1}$ is a positive constant depending on $M$ but independent of $\mu$, $\ve$, $\eta_{-}$, and $\chi_{-}$.
Furthermore, it is clear that
$
{v_+}/{2}
<
v_++\eta_m^0
<
{3v_+}/{2}
$
if $\delta\leq v_+/2M$.
Then it is also seen that
\begin{gather*}
\|\cR^{0}(\cdot,\eta^{0}_{m},\alpha^{0}(\eta_{-},\chi_{-},\eta^{0}_{m}))\| \leq C_{2} \delta^{2},
\end{gather*}
where $C_{2}$ is a positive constant depending on $M$ but independent of $\mu$, $\ve$, $\eta_{-}$, and $\chi_{-}$.
From this and \eqref{ineqG2} with $G_{*}=G^{0}$, it follows that 
\begin{gather}\label{R0-es1}
\left\| \frac{1}{\ve  \mu} \int_{1}^{\cdot}G^{0}(\cdot,s)\cR^{0}(s,\eta^{0}_{m}(s),\alpha^{0}(\eta_{-},\chi_{-},\eta^{0}_{m}))s^{n-1} \,ds \right\| \leq C_{3} \delta^{2},
\end{gather}
where $C_{3}$ is a positive constant depending on $M$ but independent of $\mu$, $\ve$, $\eta_{-}$, and $\chi_{-}$.
From these, we deduce $\eta^{0}_{m+1} \in X(M)$ if $\delta \leq \min\{v_{+}/2M, \ M_{1}/C_{1}, \  M_{1}/C_{3}\}$.
By induction, we conclude that $\{\eta^{0}_{m}\}_{m\ge 0} \subset X(M)$ if $\delta \leq \min\{v_{+}/2M, \ M_{1}/C_{1}, \ M_{1}/C_{3} \}$.

We next prove the contraction estimate
\begin{gather*}
\|\eta^{0}_{m+2}-\eta^{0}_{m+1}\| \leq \frac{1}{2} \|\eta^{0}_{m+1}-\eta^{0}_{m}\|
\end{gather*}
for suitably small $\delta$.
Once this claim is valid, it is straightforward to see that $\{\eta^{0}_{m}\}_{m\ge 0}$ is a Cauchy sequence in ${\cal {B}}^0$ and then its limit $\eta^{0} \in {\cal {B}}^0$ is the desired solution of the integral equation \eqref{single-eta2}.
Let us show this contraction estimate. 
Suppose that $\delta \leq \min\{c_{V}v_{+}/4RM, \ M_{1}/C_{1}, \ M_{1}/C_{3}\}$.  We then see using $\{\eta^{0}_{m}\}_{m\ge 0} \subset X(M)$ that
\begin{align*}
&\left| \frac{\alpha^{0}_{2}(\eta_{-},\eta^{0}_{m+1})}{\ve} - \frac{\alpha^{0}_{2}(\eta_{-},\eta^{0}_{m})}{\ve}\right| 
\\
& \quad +\|\cR^{0}(\cdot,\eta^{0}_{m+1},\alpha^{0}(\eta_{-},\chi_{-},\eta^{0}_{m+1}))-\cR^{0}(\cdot,\eta^{0}_{m},\alpha^{0}(\eta_{-},\chi_{-},\eta^{0}_{m}))\| 
\\
&\leq C_{4} \delta \|\eta^{0}_{m+1}-\eta^{0}_{m}\|,
\end{align*}
which gives
\begin{align*}
&\left\| \dfrac{R}{\ve \mu c_{V} v_+} \left(  \frac{\alpha^{0}_{2}(\eta_{-},\eta^{0}_{m+1})}{\ve} - \frac{\alpha^{0}_{2}(\eta_{-},\eta^{0}_{m})}{\ve} \right) \int_{1}^{\cdot}G^{0}(\cdot,s) s^{n-1}\,ds \right\|
\\
& + \left\| \frac{1}{\ve  \mu} \! \int_{1}^{\cdot}\! G^{0}(\cdot,s) \! \left\{\cR^{0}(s,\eta^{0}_{m+1}(s),\alpha^{0}(\eta_{-},\chi_{-},\eta^{0}_{m+1}))-\cR^{0}(s,\eta^{0}_{m}(s),\alpha^{0}(\eta_{-},\chi_{-},\eta^{0}_{m}))\right\} s^{n-1} \,ds \right\| 
\\
& \leq C_{5} \delta  \|\eta^{0}_{m+1}-\eta^{0}_{m}\|.
\end{align*}
Here $C_{5}$ is a positive constant independent of $\mu$, $\ve$, $\eta_{-}$, and $\chi_{-}$.
It follows that the desired contraction estimate holds if $\delta \leq \delta_{0}:=\min\{c_{V}v_{+}/4RM, \ M_{1}/C_{1}, \ M_{1}/C_{3}, \ 1/2C_{5}\}$. 
The proof is complete.
\end{proof}

As a consequence of Lemma \ref{thm-single-eta2}, we obtain the solvability of the limiting problem \eqref{limpro}.

\begin{lem}\label{solvablity0}
Let $\eta^{0}$ be the solution of \eqref{single-eta2} obtained in Lemma \ref{thm-single-eta2}, and define $\chi^{0}$ from $\eta^{0}$ by \eqref{form-chi1}.
Then $(\eta^{0},\chi^{0}) \in C^{2} \times C^{2}$ solves the problem \eqref{limpro} and $\|(\eta^{0},\chi^{0})\| \leq C\delta$ holds.
\end{lem}
\begin{proof}
It follows from Lemma \ref{thm-single-eta2} and \eqref{single-eta2} that
$\eta^{0} \in C^{2}$.
Combining this with the relation \eqref{form-chi1} and $\|\eta^{0}\| \leq C\delta$ yields the desired
properties of $\chi^{0}$ and completes the proof.
\end{proof}

The following lemma provides a decay estimate for the derivatives of the solution $(\eta^{0},\chi^{0})$.

\begin{lem}\label{decay0}
Let $(\eta^{0},\chi^{0})$ be the solution of \eqref{limpro} in Lemma \ref{solvablity0}.
Then there exists a constant $C$ independent of  $\ve$ and $\delta$ such that \eqref{Deta0chi0} holds, that is,
\begin{gather*}
\|(r^{n+1}\eta_r^0,r^{n+1}\chi_r^0)\|
\le
C\frac{\delta}{\ve}.
\end{gather*}
\end{lem}
\begin{proof}
Set 
\[
M:=\|r^{n+1}\eta_r^0\|.
\]
We differentiate \eqref{single-eta2} with respect to $r$,
integrate by parts in the last two integrals, and multiply the result by $r^{n+1}$  to obtain
\begin{align}
r^{n+1} \eta^{0}_{r}(r)
={}&
-\frac{\omega^{0}}{\mu\ve}
r^{2n}G^{0}(r,1)\eta_-
+\frac{R}{\mu\ve c_{V}v_+}
\frac{\alpha^{0}}{\ve}
r^{2n}G^{0}(r,1)
\notag\\
&
-2\nu(n-1)\frac{Rv_{+}}{\mu c_{V}}
r^{2n}G^{0}(r,1)
+
2\nu n(n-1)\frac{Rv_{+}}{\mu c_{V}}
r^{2n}
\int_{1}^{r}
G^{0}(r,s)s^{-n-1}\,ds
\notag\\
&
+\frac{r^{2n}}{\mu\ve}
G^{0}(r,1)
\cR^{0}(1,\eta^{0},\alpha^{0})\bigr|_{r=1}
+\frac{r^{2n}}{\mu\ve}
\int_{1}^{r}
G^{0}(r,s)
\frac{d}{ds}
\cR^{0}(s,\eta^{0}(s),\alpha^{0})
\,ds.
\label{fDeta0}
\end{align}
In what follows, $C>0$ is a constant independent of $\delta$ and $\ve$.
We observe from \eqref{hyoka1}, \eqref{hyoka2}, and $|\cR^{0}(1,\eta^{0},\alpha^{0})\bigr|_{r=1}| \leq C\delta$ that the first, second, third, and fifth terms on the right-hand side are bounded above by $C\delta\ve^{-1}$. The fourth term is also bounded above by $C\delta$ as a consequence of \eqref{ineqG} with $G_{*}=G^{0}$ and $l=n+1$. 
What is left is to estimate the last term.
It follows that
\[
\left|
\frac{d}{ds}
\cR^0(s,\eta^0(s),\alpha^0)
\right|
\leq
C\delta(M+1)s^{-n-1}.
\]
Hence, by \eqref{ineqG} with \(G_*=G^0\) and $l=n+1$,
\begin{align*}
\frac{r^{2n}}{\mu\ve}
\left|
\int_1^r
G^0(r,s)
\frac{d}{ds}
\cR^0(s,\eta^0(s),\alpha^0)\,ds
\right|
\leq
\frac{C\delta(M+1)}{\ve}
r^{2n}
\int_1^rG^0(r,s)s^{-n-1}\,ds
\leq
C\delta(M+1).
\end{align*}
Substituting the preceding estimates into \eqref{fDeta0}, we obtain
$
M
\leq
C\delta\ve^{-1}
+
C\delta(M+1).
$
Taking \(\delta\) sufficiently small, we absorb the last term
into the left-hand side and obtain
$
\|r^{n+1}\eta_r^0\|
\leq
C\delta\ve^{-1}.
$
Finally, differentiating \eqref{form-chi1} and using the preceding
estimate, 
we obtain
$
\|r^{n+1}\chi_r^0\|
\leq
C\delta\ve^{-1}.
$
This proves \eqref{Deta0chi0}.
\end{proof}


We next investigate the far-field behavior of the solution to the limiting problem.

\begin{lem}\label{cor:eta-limit}
Let $(\eta^{0},\chi^{0})$ be the solution of \eqref{limpro} obtained in Lemma \ref{solvablity0}.
Then the far-field limits $\eta_{\infty}^{0}=\lim_{r\to\infty}\eta^{0}(r)$ and 
$\chi_{\infty}^{0}=\lim_{r\to\infty}\chi^{0}(r)$ exist, and the relation \eqref{farfield-relation2} holds, that is,
\begin{gather*}
\theta_+\eta_{\infty}^{0}
=
v_+\chi_{\infty}^{0}.
\end{gather*}
\end{lem}

\begin{proof}
Due to \eqref{Deta0chi0} and $r^{-(n+1)} \in L^{1}(1,\infty)$, both $\eta^0(r)$ and
$\chi^0(r)$ have finite limits as $r\to\infty$. 
We next take the limit $r\to\infty$ in the first equation of \eqref{limpro} using the fact ${\eta_r^0(r)}/{r^{n-1}} \to 0$ to obtain
\[
0=
\frac{R}{v_+}\chi^{0}_\infty
-\frac{R\chi^{0}_\infty\eta^{0}_\infty}
{v_+(v_++\eta^{0}_\infty)}
-\frac{R\theta_+}{v_+^2}\eta^{0}_\infty
+\frac{R\theta_+(\eta_\infty^{0})^2}
{v_+^2(v_++\eta^{0}_\infty)}.
\]
From this and $v_++\eta^{0}_\infty>0$, we arrive at \eqref{farfield-relation2}.
\end{proof}

We now complete the proof of Proposition \ref{pro3.2}.
\begin{proof}[Proof of Proposition \ref{pro3.2}.]
Proposition \ref{pro3.2} immediately follows from Lemmas \ref{thm-single-eta2}--\ref{cor:eta-limit}. 
\end{proof}

\section{Far field layer}\label{S6}
In this section, we construct the far-field layer.
To identify the asymptotic structure in the far-field region,
we introduce the stretched variable
$\bar r=\kappa^{1/(n-2)}r$.
Applying this scaling to the first two equations of
\eqref{reforme} and formally letting $\kappa\to0$, we obtain
\begin{align*}
 \left\{
\begin{aligned}
& \bar{\eta}^{\rm far} = \dfrac{v_+}{\theta_+}\bar{\chi}^{\rm far},
\\
&\bar{\chi}_{\bar{r}}^{\rm far} 
-\frac{\ve c_V}{\bar{r}^{n-1}}\bar{\chi}^{\rm far}
= 
-\frac{\bar{\alpha}^{\rm far}}{\bar{r}^{n-1}} 
+ \dfrac{\ve R\theta_{+}}{v_{+} \bar{r}^{n-1}}\bar{\eta}^{\rm far}.
 \end{aligned}
 \right.
\end{align*}
The first equation shows that the pressure remains constant to the first order in the stretched variables.
Substituting this relation into the second equation yields a closed equation for the far-field profile.

Motivated by this limiting system, we define the far-field layer as the solution of
\begin{align*}
 \left\{
\begin{aligned}
&\bar{\chi}_{\bar r}^{\mathrm{far}}
-\frac{\varepsilon(c_V+R)}{\bar r^{\,n-1}}
\bar{\chi}^{\mathrm{far}}
=
-\frac{\bar{\alpha}^{\mathrm{far}}}{\bar r^{\,n-1}},
\\
&\bar{\chi}^{\mathrm{far}}(\kappa^{1/(n-2)})=\chi^{0}_\infty,
\\
&\lim_{\bar r\to\infty}\bar{\chi}^{\mathrm{far}}(\bar r)
=0,
\end{aligned}
\right.
\end{align*}
where $\chi^{0}_\infty$ is defined in \eqref{eta-chi-infty1}.
The corresponding far-field layer in the original variables is defined by
\begin{gather*}
\tilde{\chi}^{\mathrm{far}}(r)
:=
\bar{\chi}^{\mathrm{far}}(\kappa^{1/(n-2)}r),
\quad
\tilde{\alpha}^{\mathrm{far}}
:=
\bar{\alpha}^{\mathrm{far}}.
\end{gather*}
These functions admit the explicit representations
\begin{align*}
\tilde{\chi}^{\mathrm{far}}(r)
&=
\chi^{0}_\infty
\left(
1-e^{-\frac{\varepsilon\sigma^{\mathrm{far}}}
{\kappa(n-2)}}
\right)^{-1}
\left(
1-e^{-\frac{\varepsilon\sigma^{\mathrm{far}}}
{\kappa(n-2)r^{n-2}}}
\right),
\\
\tilde{\alpha}^{\mathrm{far}}
&=
\varepsilon(c_V+R)\chi^{0}_\infty
\left(
1-e^{-\frac{\varepsilon\sigma^{\mathrm{far}}}
{\kappa(n-2)}}
\right)^{-1},
\end{align*}
where
\begin{gather*}
\sigma^{\mathrm{far}}
:=c_V+R.
\end{gather*}

For the subsequent analysis, it is convenient to subtract the limiting state $(\eta^{0}_{\infty}, \chi^{0}_{\infty})$ defined in \eqref{eta-chi-infty1} and redefine the far-field layer by
\begin{align*}
\chi^{\rm far}
&:=
\tilde{\chi}^{\rm far}-\chi^{0}_\infty=
-\chi^{0}_\infty
\left(
1-e^{-\frac{\varepsilon\sigma^{\rm far}}
{\kappa(n-2)}}
\right)^{-1}
\left(
e^{-\frac{\varepsilon\sigma^{\rm far}}
{\kappa(n-2)r^{n-2}}}
-
e^{-\frac{\varepsilon\sigma^{\rm far}}
{\kappa(n-2)}}
\right),
\\
\alpha^{\rm far}
&:=
\tilde{\alpha}^{\rm far}
-\varepsilon(c_V+R)\chi^{0}_\infty
=
\varepsilon(c_V+R)\chi^{0}_\infty
e^{-\frac{\varepsilon\sigma^{\rm far}}
{\kappa(n-2)}}
\left(1-
e^{-\frac{\varepsilon\sigma^{\rm far}}
{\kappa(n-2)}}
\right)^{-1},
\\
\eta^{\rm far}
&:=
\tilde{\eta}^{\rm far}-\eta^{0}_\infty
=
\frac{v_+}{\theta_+}\chi^{\rm far},
\end{align*}
where the last equality follows from
\eqref{farfield-relation2}.
Note that $(\eta^{\rm far},\chi^{\rm far})$ coincide with the pair defined in \eqref{FFL1},
and the triple
$(\eta^{\rm far},\chi^{\rm far},\alpha^{\rm far})$
satisfies the following system:
\begin{align}\label{FFLsys1}
 \left\{
\begin{aligned}
&\frac{R\theta_+r^{n-1}}
{\ve\mu v_+^2}
\eta^{\rm far}
=
\frac{Rr^{n-1}}
{\ve\mu v_+}
\chi^{\rm far},
\\
&(\chi^{\rm far})_{r}
-\frac{\ve c_V}
{\kappa r^{n-1}}
\chi^{\rm far}
=
-\frac{\alpha^{\rm far}}
{\kappa r^{n-1}}
+\frac{\ve R\theta_+}
{\kappa v_+r^{n-1}}
\eta^{\rm far},
\quad r>1,
\\
&\eta^{\rm far}(1)=0,
\qquad
\chi^{\rm far}(1)=0,
\\
&\lim_{r\to\infty}
\eta^{\rm far}(r)
=
-\eta^{0}_\infty,
\qquad
\lim_{r\to\infty}
\chi^{\rm far}(r)
=
-\chi^{0}_\infty.
\end{aligned}
\right.
\end{align}

We derive several estimates that will be used in the analysis of the vanishing heat conductivity limit.

\begin{lem}\label{lem:far-weighted-decay}
Let \(n\geq 3\) and \(q>0\). If $\kappa \leq \ve\sigma^{\rm far} $,
then
\begin{gather}
\label{eta-far-weighted}
\sup_{r\in [1,\infty)}r^{-q}|\eta^{\rm far}(r)|
\leq
C\delta
\left(
\frac{\kappa}{\ve}
\right)^{\frac{q}{n-2}},
\\
\label{eta-far-integral}
\sup_{r\in [1,\infty)}
r^{\,n-2}
\int_r^\infty
\frac{|\eta^{\rm far}(s)|}{s^{2n-1}}\,ds
\le
C\delta
\left(
\frac{\kappa}{\ve}
\right)^{\frac{n}{n-2}},
\end{gather}
where \(C>0\) is a constant independent of \(\kappa\), \(\ve\), and \(\delta\).
\end{lem}

\begin{proof}
Set
\[
a
:=
\frac{\ve\sigma^{\rm far}}{(n-2)\kappa}.
\]
By assumption, \(a \geq 1/(n-2)\) holds, and hence we have
$
1-e^{-a}
\geq
1-e^{-\frac{1}{n-2}}>0.
$
From this and the definition of $\eta^{\rm far}$, it follows that
\begin{align}
r^{-q}|\eta^{\rm far}(r)|
&\leq
C|\eta^{0}_\infty|r^{-q}
\left(
e^{-\frac{a}{r^{n-2}}}
+
e^{-a}
\right).
\label{eta-far-weighted-1}
\end{align}
We first estimate the first term on the right hand side.
Introducing the variable
$
x=a r^{2-n}
$,
we obtain
$r^{-q}e^{-\frac{a}{r^{n-2}}}=a^{-\frac{q}{n-2}} x^{\frac{q}{n-2}}e^{-x}$.
The function $x^{\frac{q}{n-2}}e^{-x}$ attains its maximum at $x=q/(n-2)$. 
For the second term in \eqref{eta-far-weighted-1}, it follows from Taylor's theorem that 
$
r^{-q}e^{-a}
\leq
e^{-a}
\leq
Ca^{-\frac{q}{n-2}}.
$
Recall that $|\eta^{0}_{\infty}| \leq \|\eta^{0}\| \leq C\delta$. 
Combining these, we conclude \eqref{eta-far-weighted}. 

Finally, choosing \(q=n\) in \eqref{eta-far-weighted}, we obtain
\[
|\eta^{\rm far}(s)|=
s^{n}s^{-n}|\eta^{\rm far}(s)|
\le
C\delta
\left(
\frac{\kappa}{\ve}
\right)^{\frac{n}{n-2}}
s^n.
\]
Therefore,
\begin{align*}
\int_r^\infty
\frac{|\eta^{\rm far}(s)|}{s^{2n-1}}\,ds
&\le
C\delta
\left(
\frac{\kappa}{\ve}
\right)^{\frac{n}{n-2}}
r^{-(n-2)}.
\end{align*}
This proves \eqref{eta-far-integral}.
\end{proof}

\section{Vanishing Heat Conductivity Limit}\label{S7}
In this section, we study the vanishing heat conductivity limit $\kappa \to 0$.
We set
\begin{gather*}
\bar{\eta}:=\eta-\eta^{0}-\eta^{\rm far}, \quad
\bar{\chi}:=\chi-\chi^{0}-\chi^{\rm far}, \quad
\bar{\alpha}:=\alpha-\alpha^{0}-\alpha^{\rm far},
\end{gather*}
and then see from \eqref{reforme}, \eqref{limpro}, and \eqref{FFLsys1} that the triple $(\bar{\eta},\bar{\chi},\bar{\alpha})$ satisfies
\begin{align}
\label{bar-system}
\left\{
\begin{aligned}
&\bar{\eta}_{r}
+\frac{R\theta_+r^{n-1}}{\ve\mu v_+^2}\bar{\eta}
=
\frac{Rr^{n-1}}{\ve\mu v_+}\bar{\chi}
-\eta^{\rm far}_{r}
+\frac{r^{n-1}}{\ve\mu}\mathcal{N}^{\eta},
\\
&\bar{\chi}_{r}
-\frac{\ve c_V}{\kappa r^{n-1}}\bar{\chi}
=
-\frac{\bar{\alpha}}{\kappa r^{n-1}}
+\frac{\ve R\theta_+}{\kappa v_+r^{n-1}}\bar{\eta}
-\chi_r^0
+\frac{\ve}{\kappa r^{n-1}}\mathcal{N}^{\chi},
\quad r>1,
\\
&\bar{\eta}(1)=0,\quad
\bar{\chi}(1)=0,
\\
&\lim_{r\to\infty}\bar{\eta}(r)=0,\qquad
\lim_{r\to\infty}\bar{\chi}(r)=0.
\end{aligned}
\right.
\end{align}
Here,
\begin{align*}
\mathcal N^\eta
&:=
-\frac{R\chi\eta}{v_+(v_++\eta)}
+\frac{R\chi^0\eta^0}{v_+(v_++\eta^0)}
+F(r,\eta)-F(r,\eta^0),
\\
\mathcal N^\chi
&:=H(r,\eta)-H(r,\eta^0).
\end{align*}
Using the kernels \(G\) and \(K\), defined in  \eqref{G1} and \eqref{K1}, and the conditions 
$\bar{\eta}(1)=\lim_{r\to\infty}\bar{\chi}(r)=0$, we can write $\bar{\eta}$ and $\bar{\chi}$ as 
\begin{align}
\bar{\eta}(r)
&=
\int_{1}^{r}G(r,s)
\left\{
\frac{R}{\varepsilon\mu v_{+}}
 s^{n-1}\bar{\chi}(s)
-\eta^{\rm far}_{s}(s)
+\frac{s^{n-1}}{\varepsilon\mu}
 \mathcal N^{\eta}(s)
\right\}\,ds,
\label{bar-eta-formula}
\\
\bar{\chi}(r)
&=
\int_{r}^{\infty}K(r,s)
\left\{
\frac{\bar{\alpha}}{\kappa s^{n-1}}
-\frac{\varepsilon R\theta_{+}}
 {\kappa v_{+}s^{n-1}}\bar{\eta}(s)
+\chi^{0}_{s}(s)
-\frac{\varepsilon}{\kappa s^{n-1}}
 \mathcal N^{\chi}(s)
\right\}\,ds.
\label{bar-chi-formula}
\end{align}
Furthermore, the formula \eqref{bar-chi-formula} with the condition $\bar{\chi}(1)=0$ gives
\begin{align}
\frac{\bar\alpha}{\varepsilon c_V}
&=
\left\{
\frac{\varepsilon R\theta_+}{\kappa v_+}
\int_1^\infty
K(1,s)\frac{\bar\eta(s)}{s^{n-1}}\,ds
-\int_1^\infty
K(1,s)\chi_s^0(s)\,ds
+\frac{\varepsilon}{\kappa}
\int_1^\infty
K(1,s)
\frac{\mathcal N^\chi(s)}{s^{n-1}}\,ds
\right\}
\notag \\
& \quad
\times \left(
1-
e^{-\frac{c_V\varepsilon}{\kappa(n-2)}}
\right)^{-1}.
\label{bar-alpha-formula}
\end{align}

We begin with estimating the terms $\mathcal N^\eta$ and $\mathcal N^\chi$.

\begin{lem}\label{est-N}
Let $n \ge 3$ and $\gamma<2$. For any $v_+>0$ and $\theta_+>0$, there exists  a positive constant $\delta_0$, independent of $\ve$, $\delta$, and $\kappa$, such that if $\delta + \kappa/\ve^{3} \leq \delta_0$, then
\begin{align}
|\mathcal N^\eta| & \leq C_{2}\delta
\left(
\|(\bar\chi,\bar\eta)\|
+\frac{|\eta^{\rm far}|}{\ve r^n}
+\int_r^\infty
\frac{|\eta^{\rm far}(s)|}{s^{2n-1}}\,ds
\right),
\label{est-Neta} \\
|\mathcal N^\chi| & \leq C_{2}\delta \left(\|\bar{\eta}\| + \frac{|\eta^{\rm far}|}{r^{n}} +\int_r^\infty
\frac{|\eta^{\rm far}(s)|}{s^{2n-1}}\,ds\right).
\label{est-Nchi}
\end{align}
\end{lem}
\begin{proof}
We only prove \eqref{est-Neta}, since it is easy to prove \eqref{est-Nchi}.
To show \eqref{est-Neta}, we write $\mathcal N^\eta$ as
\begin{align*}
\mathcal N^\eta
&=
-\frac{R}{v_+}
\left(
\frac{\chi\eta}{v_++\eta}
-
\frac{\chi^0\eta^0}{v_++\eta^0}
\right)
+
\frac{R\theta_+}{v_+^2}
\left(
\frac{\eta^2}{v_++\eta}
-
\frac{(\eta^0)^2}{v_++\eta^0}
\right)
\\
&\quad +
\frac{\varepsilon^2}{r^{2(n-1)}}
\bigl(\bar\eta+\eta^{\rm far}\bigr)
-
\varepsilon^2(n-1)
\int_r^\infty
\frac{\bar\eta(s)+\eta^{\rm far}(s)}{s^{2n-1}}\,ds
\\
&
=:I_{1}+I_{2}+I_{3}+I_{4}.
\end{align*}
Recall the relations in \eqref{farfield-relation2} and \eqref{FFLsys1}, i.e.
\[
\chi^{0}_{\infty}=\frac{\theta_+}{v_+}\eta^{0}_{\infty}, \quad
\chi^{\rm far}=\frac{\theta_+}{v_+}\eta^{\rm far}.
\]
The first two terms $I_{1}$ and $I_{2}$ can be rewritten by the definitions of $\bar\eta$ and $\bar\chi$ and the relations as
\begin{align*}
I_{1}+I_{2}&=-\frac{R}{v_+}
\left[
\frac{\chi\eta}{v_++\eta}
-
\frac{\chi^0\eta^0}{v_++\eta^0}
-
\frac{\theta_+}{v_+}
\left(
\frac{\eta^2}{v_++\eta}
-
\frac{(\eta^0)^2}{v_++\eta^0}
\right)
\right]
\\
&=
-\frac{R}{v_+}
\left[
\frac{\left(\chi-\frac{\theta_+}{v_+}\eta\right)\eta}{v_++\eta}
-
\frac{\left(\chi^0-\frac{\theta_+}{v_+}\eta^0\right)\eta^0}{v_++\eta^0}
\right]
\\
&=
-\frac{R}{v_+}
\left[
\left(\bar\chi-\frac{\theta_+}{v_+}\bar\eta\right)
\frac{\eta}{v_++\eta}
+
\left(\chi^0-\frac{\theta_+}{v_+}\eta^0\right)
\frac{v_+\bigl(\bar\eta+\eta^{\rm far}\bigr)}
{(v_++\eta)(v_++\eta^0)}
\right].
\end{align*}
We observe using the relation and \eqref{Deta0chi0} that 
\begin{align*}
\left|
\chi^0-\frac{\theta_+}{v_+}\eta^0
\right|
=
\left|
\chi^0-\chi^{0}_{\infty}-\frac{\theta_+}{v_+}(\eta^0-\eta^{0}_{\infty})
\right|
=
\left|
\int^{\infty}_{r} \chi^{0}_{s} ds -
\frac{\theta_+}{v_+}\int^{\infty}_{r} \eta^{0}_{s} ds
\right|
\leq
C\frac{\delta}{\ve r^{n}}.
\end{align*}
From this and 
$|(\eta,\chi,\eta^{0},\chi^{0})|\leq C\delta$, we arrive at
\begin{align*}
&\left|
I_{1}+I_{2}
\right|
\leq
C\delta
\left(
|\bar\chi|+|\bar\eta|
+
\frac{|\eta^{\rm far}|}{\ve r^n}
\right).
\end{align*}
It is straightforward to estimate $|I_{3}|$ by the right hand side of \eqref{est-Neta}.
Finally, the last term \(I_{4}\) can be estimated as
\begin{align*}
|I_{4}|
&\leq
\varepsilon^2(n-1)
\int_r^\infty
\frac{|\bar\eta(s)|}{s^{2n-1}}\,ds
+
\varepsilon^2(n-1)
\int_r^\infty
\frac{|\eta^{\rm far}(s)|}{s^{2n-1}}\,ds
\\
&\leq
C\delta
\left(
\|\bar\eta\|
+
\int_r^\infty
\frac{|\eta^{\rm far}(s)|}{s^{2n-1}}\,ds
\right).
\end{align*}
Combining the above estimates, we conclude \eqref{est-Neta}.
\end{proof}

We are now in a position to prove Theorem \ref{thm3.3}.

\begin{proof}[Proof of Theorem \ref{thm3.3}]
We first derive an estimate for \(\bar\chi\) in the same way as in
\eqref{l-chi-es}.  Applying \eqref{hyokaK2} with \(K_*=K\) to \eqref{bar-chi-formula},
we obtain
\begin{align}
\|\bar\chi\|
&\leq
\left|\frac{\bar\alpha}{\varepsilon c_V}\right|
+\frac{R\theta_+}{c_Vv_+}\|\bar\eta\|
+C
 \left\|\int_{\cdot}^\infty K(\cdot,s)\chi_s^0(s)\,ds\right\|
+C
 \left\|
 \frac{\varepsilon}
 {\kappa}
 \int_{\cdot}^\infty
 K(\cdot,s)\frac{\mathcal N^\chi}
 {s^{n-1}}\,ds
 \right\|.
\label{bar-chi-pre}
\end{align}
From \eqref{Deta0chi0} and \eqref{hyokaK2} with \(K_*=K\), it follows that
\begin{equation}\label{chi0-source-est}
\left\|\int_{\cdot}^\infty K(\cdot,s)\chi_s^0(s)\,ds\right\|
\leq
C\delta\frac{\kappa}{\varepsilon^2}.
\end{equation}
Moreover, using Lemmas \ref{est-N} and \ref{lem:far-weighted-decay}, and then applying \eqref{hyokaK2} with \(K_*=K\), we have
\begin{align}
 \left\|
 \frac{\varepsilon}
 {\kappa}
 \int_{\cdot}^\infty
 K(\cdot,s)\frac{\mathcal N^\chi}
 {s^{n-1}}\,ds
 \right\|
\leq 
C\delta\|(\bar\eta,\bar\chi)\|
+C\delta^2
 \left(\frac{\kappa}{\varepsilon}\right)^{\frac{n}{n-2}}.
\label{Nchi-source-est}
\end{align}
Combining
\eqref{bar-chi-pre}--\eqref{Nchi-source-est},
we conclude that
\begin{align}
\|\bar\chi\|
&\le
\left|\frac{\bar\alpha}{\varepsilon c_V}\right|
+\frac{R\theta_+}{c_Vv_+}\|\bar\eta\|
+
C\delta\|(\bar\eta,\bar\chi)\|
+
C\delta\frac{\kappa}{\varepsilon^{2}}
+
C\delta^2
\left(\frac{\kappa}{\varepsilon}\right)^{\frac{n}{n-2}}.
\label{bar-eta-from-chi}
\end{align}

Next, we estimate $\bar\eta$. Applying \eqref{ineqG2} with $G_*=G$ to
\eqref{bar-eta-formula}, we obtain
\begin{align}
\|\bar\eta\|
&\le
\frac{v_+}{\theta_+}\|\bar\chi\|
+
\left\|
\int_1^{\cdot}
G(\cdot,s)\eta_s^{\rm far}(s)\,ds
\right\|
+
\left\|
\frac1{\varepsilon\mu}
\int_1^{\cdot}
G(\cdot,s)s^{n-1}\mathcal N^\eta(s)\,ds
\right\| .
\label{bar-eta-pre}
\end{align}
Note that 
\[
\eta^{\rm far}_r(r)=-C^{\rm far}\partial_rK^{\rm far}(r,\infty), 
\]
where
\[
K^{\mathrm{far}}(r,s)
:=
e^{
-\frac{\varepsilon\sigma^{\mathrm{far}}}
{\kappa(n-2)}
\left(
\frac1{r^{n-2}}
-
\frac1{s^{n-2}}
\right)},
\quad 
C^{\rm far}
:=
\frac{v_+}{\theta_+}\chi_\infty
\left(1-
e^{-\frac{\varepsilon\sigma^{\rm far}}
{\kappa(n-2)}}
\right)^{-1}.
\]
By Lemma~\ref{lem:GK-mixed} with $G_{*}=G$ and $K_{*}=K^{\rm far}$,
\begin{align}
\left\|
\int_1^{\cdot}
G(\cdot,s)\eta_s^{\rm far}(s)\,ds
\right\|
\le
C\delta\varepsilon
\left(\frac{\kappa}{\varepsilon}\right)^{\frac{n}{n-2}}
\le
C\delta^{2}
\left(\frac{\kappa}{\varepsilon}\right)^{\frac{n}{n-2}}.
\label{etafar-source-est}
\end{align}
Moreover, using Lemmas \ref{est-N} and \ref{lem:far-weighted-decay}, and then applying \eqref{ineqG2} with $G_{*}=G$, we obtain
\begin{align}
\left\|
\frac1{\varepsilon\mu}
\int_1^{\cdot}
G(\cdot,s)s^{n-1}\mathcal N^\eta\,ds
\right\|
&\le
C\delta\|(\bar\eta,\bar\chi)\|
+
C\delta^2
\left(\frac{\kappa}{\varepsilon}\right)^{\frac{n}{n-2}}
+
C\delta^2
\frac1{\varepsilon}
\left(\frac{\kappa}{\varepsilon}\right)^{\frac{n}{n-2}}.
\label{Neta-source-est}
\end{align}
Combining
\eqref{bar-eta-pre}--\eqref{Neta-source-est},
we conclude that
\begin{align}
\|\bar\eta\|
&\le
\frac{v_+}{\theta_+}\|\bar\chi\|
+
C\delta\|(\bar\eta,\bar\chi)\|
+
C\delta^2
\left(\frac{\kappa}{\varepsilon}\right)^{\frac{n}{n-2}}
+
C\delta^2
\frac1{\varepsilon}
\left(\frac{\kappa}{\varepsilon}\right)^{\frac{n}{n-2}}.
\label{bar-chi-from-eta}
\end{align}

Substituting \eqref{bar-chi-from-eta} into \eqref{bar-eta-from-chi}, and
using \(R/c_V=\gamma-1\), we arrive at
\begin{align}
(2-\gamma)\|\bar\chi\|
&\leq
\left|\frac{\bar\alpha}{\varepsilon c_V}\right|
+C\delta\|(\bar\eta,\bar\chi)\|
+C\delta\frac{\kappa}{\varepsilon^2}
+C\delta^2
 \left(\frac{\kappa}{\varepsilon}\right)^{\frac{n}{n-2}}
+C\delta^2\frac1{\varepsilon}
 \left(\frac{\kappa}{\varepsilon}\right)^{\frac{n}{n-2}}.
\label{bar-chi-est}
\end{align}
This is the analogue of \eqref{l-chi-es} for the remainder
\(\bar\chi\).
We next substitute \eqref{bar-chi-est} into
\eqref{bar-chi-from-eta} to obtain
\begin{align}
(2-\gamma)\|\bar\eta\|
&\leq
\frac{v_+}{\theta_+}
\left|\frac{\bar\alpha}{\varepsilon c_V}\right|
+C\delta\|(\bar\eta,\bar\chi)\|
+C\delta\frac{\kappa}{\varepsilon^2}
+C\delta^2
 \left(\frac{\kappa}{\varepsilon}\right)^{\frac{n}{n-2}}
+C\delta^2\frac1{\varepsilon}
 \left(\frac{\kappa}{\varepsilon}\right)^{\frac{n}{n-2}}.
\label{bar-eta-est}
\end{align}
This is the counterpart of \eqref{l-eta-es}.
It remains to estimate \(\bar\alpha\).  Recalling the formula \eqref{bar-alpha-formula} and repeating the proof of
\eqref{l-alpha-es}, we apply \eqref{ker2} with the aid of the condition $\bar{\eta}(1)=0$ and the first equation in
\eqref{bar-system} and also use \eqref{chi0-source-est}--\eqref{Nchi-source-est} to obtain
\begin{align}
\left|\frac{\bar\alpha}{\varepsilon c_V}\right|
&\leq
\tilde{C}_0
\left(1+\frac1{\varepsilon}\right)
\left(\frac{\kappa}{\varepsilon}\right)^{1/2}
\|(\bar\eta,\bar\chi)\|
+C\delta\|(\bar\eta,\bar\chi)\|
\nonumber\\
&\quad
+C\delta\frac{\kappa}{\varepsilon^2}
+C\delta^2
 \left(\frac{\kappa}{\varepsilon}\right)^{\frac{n}{n-2}}
+C\delta^2\frac1{\varepsilon}
 \left(\frac{\kappa}{\varepsilon}\right)^{\frac{n}{n-2}},
\label{bar-alpha-est}
\end{align}
where $\tilde{C}_{0}>0$ is a constant independent of $\ve$, $\delta$, and $\kappa$.

We now simplify the exponent of $\kappa$. 
Since \(n/(n-2)\geq1\), \(\kappa/\varepsilon\leq1\), and \(\kappa/\varepsilon^{3}\leq1\), we have
\begin{align*}
\delta \frac{\kappa}{\varepsilon^2} \leq \delta \kappa^{1/3}, \quad
\delta^2
\left(\frac{\kappa}{\varepsilon}\right)^{\frac{n}{n-2}}
\leq
\delta^2\frac{\kappa}{\varepsilon}
\leq
\delta^2\kappa^{1/3},
\quad
\frac{\delta^2}{\varepsilon}
\left(\frac{\kappa}{\varepsilon}\right)^{\frac{n}{n-2}}
\leq
\delta^2\frac{\kappa}{\varepsilon^2}
\leq
\delta^2\kappa^{1/3}.
\end{align*}
Therefore, combining \eqref{bar-chi-est}--\eqref{bar-alpha-est} and taking \(\delta \ll1 \) and $\kappa/\ve^{3} \ll 1$, we arrive at
\begin{equation*}
\|(\bar\eta,\bar\chi)\|
+
\left|\frac{\bar\alpha}{\varepsilon c_V}\right|
\leq
C\delta \kappa^{1/3}.
\end{equation*}
This completes the proof.
\end{proof}

\noindent
{\bf Funding.}
The work of M. Suzuki was supported by JSPS KAKENHI Number 25K00915.

\medskip
\noindent
{\bf Conflict of Interest.}
The authors declare that they have no conflict of interest.

\medskip
\noindent
{\bf Data Availability.}
Data sharing is not applicable to this article as no datasets were generated or analyzed during the current study.


\end{document}